\documentclass[12pt]{article}

\usepackage{latexsym}
\usepackage{bm}
\usepackage{amsfonts}
\usepackage{amssymb}
\usepackage{amsmath}
\usepackage{xspace}

\usepackage[cm]{fullpage}

\usepackage{amsmath,amscd, float,rotating}
\usepackage{amscd}

\usepackage{color}

\usepackage{latexsym}
\usepackage{amsfonts}
\usepackage{amsmath}
\usepackage{amssymb}

\numberwithin{equation}{section}

\def\XXint#1#2#3{{\setbox0=\hbox{$#1{#2#3}{\int}$}
\vcenter{\hbox{$#2#3$}}\kern-.5\wd0}}

\newcommand{\h}[2]{h^{#1 \overline{#2}}}

\def\endpf{\hbox{\vrule height1.5ex width.5em}}
\def\h{\hbox}
\def\a{\alpha}
\def\b{\beta}

\newcommand{\sm}{\setminus}

\newcommand{\CC}{{\mathbb C}}
\newcommand{\RR}{{\mathbb R}}

\newcommand{\BB}{{\mathbb B}}

\newcommand{\ra}{\rightarrow}

\def \-{\overline}
\def\wh{\widehat}
\def\wt{\widetilde}

\newcommand{\p}{\partial}

\def\endpf{\hbox{\vrule height1.5ex width.5em}}

\begin{document}
\newcounter{remark}
\newcounter{theor}
\setcounter{remark}{0} \setcounter{theor}{1}
\newtheorem{claim}{Claim}
\newtheorem{theorem}{Theorem}[section]
\newtheorem{proposition}{Proposition}[section]
\newtheorem{lemma}{Lemma}[section]
\newtheorem{question}{Question}[section]
\newtheorem{definition}{Definition}[section]
\newtheorem{corollary}{Corollary}[section]
\newtheorem{example}{Example}[section]
\newenvironment{proof}[1][Proof]{\begin{trivlist}
\item[\hskip \labelsep {\bfseries #1}]}{\end{trivlist}}
\newenvironment{remark}[1][Remark]{\addtocounter{remark}{1} \begin{trivlist}
\item[\hskip \labelsep {\bfseries #1
\thesection.\theremark}]}{\end{trivlist}}


\centerline{\bf \Large Holomorphic isometry from a K\"ahler
manifold} \centerline{\bf \Large   into a product of complex
projective manifolds}


\medskip
\begin{center}{\large Xiaojun Huang$^*$\footnote{Supported in part by  NSF-1101481.} ~  and  ~  Yuan Yuan$^\dagger$ }
\end{center}


\medskip

\begin{abstract} {\small We  study the global property of local holomorphic
isometric mappings from a class of K\"ahler manifolds into a product
of projective algebraic manifolds with induced Fubini-Study metrics,
where isometric  factors are allowed to be negative.}
\end{abstract}

\bigskip\bigskip
\section{Introduction}


In this paper, we study the global property for local holomorphic
isometric mappings, up to isometric factors that are allowed  even
to be negative, from a K\"ahler manifold
into a Cartesian product of projective algebraic manifolds equipped
with the induced Fubini-Study metrics.
Investigations of this kind started with a  paper of Calabi, who
first studied the global extension and Bonnet type rigidity of a
local holomorphic isometric embedding from  complex manifolds with
real analytic K\"{a}hler metrics \cite{C} into complex space forms.
Afterwards, there appeared many studies  along these lines of
research.
(see \cite{U} [DL], for instance). In 2003, motivated from problems
in algebraic number theory, Clozel-Ullmo \cite{CU} were led to
consider rigidity problems for local holomorphic isometries between
bounded symmetric domains equipped with their Bergman metrics. More
precisely, by reducing the problem to the rigidity problem for local
holomorphic isometries, they proved an algebraic correspondence in
the quotient of a bounded symmetric domain preserving the Bergman
metric has to be a modular correspondence in the case of unit disc
in the complex plane and bounded symmetric domains of rank $\geq 2$.
More recently, Mok carried out a  systematic study of this type of
problems in a very general setting. Many  important results have
been obtained by Mok and Mok-Ng. (See [M2] [N1-2], in particular,
the papers by Mok [M3] and Mok-Ng \cite{MN}, and the references
therein). Mok in [M2-3] proved the total geodesy for a local
holomorphic isometric embedding between bounded symmetric domains
$D$ and $\Omega$ when either (i) the rank of each irreducible
component of $D$ is at least two  or (ii) $D=\BB^n$ and
$\Omega=(\BB^n)^p$ for $n\geq 2$. Mok and Ng in [MN] proved the
total geodesy  when the map is a local volume preserving map, which
in particular has applications to answer, in the affirmative,
questions of Clozel-Ullmo in algebraic number theory. In a recent
joint paper of the second author with Zhang [YZ], the total geodesy
is  obtained in the case of $D=\BB^n$ and $\Omega=\BB^{N_1}
\times\cdots \times \BB^{N_p}$ with $n\geq 2$ and $N_l$ arbitrary
for $1 \leq l \leq p$. Earlier, Ng in [N2] had established similar
result when $p=2$ and $2 \le n\le N_1, N_2\le 2n-1.$


\medskip

When manifolds are  Hermitian symmetric spaces of compact type,
it is  well known that the total geodesy   for (local) holomorphic
 isometric embeddings  is no longer the
case, even for mappings between complex projective spaces equipped
with the standard Fubini-Study metrics.  For instance, the Veronese
embedding from $(\mathbb{P}^n, \omega_{n})$ into
$(\mathbb{P}^{\frac{n(n+3)}{2}}, \omega_{\frac{n(n+3)}{2}})$ is an
isometric embedding with conformal constant 2, which is not linear.
In this setting, the rigidity problem is to find out when the maps
are holomorphic isometries up to isometric constants.

In this paper, we carry out  a study  of the rigidity problem for
local holomorphic conformal maps into the product of complex
projective manifolds with the induced Fubini-Study metrics.
Our  conformal factors are allowed to
  have mixed signs.
 This has immediate application to the case when the target manifold is the product of  the Hermitian symmetric
 spaces of compact type equipped with the canonical metrics.

Geometrically, such a problem may be regarded as the question of
understanding what modification of the canonical metric on the
source manifold admits only rigid holomorphic isometric embeddings
(up to conformal factors) into product of projective spaces. To
state our main result,
we let
$$\ \  (M,\omega_m),\ (M_1,\omega_{M_1}),\ \ \cdots,\ \
(M_m,\omega_{M_m}),\ \ (M_1',\omega_{M'_1}),\cdots,\ \
(M_v',\omega_{M'_v})$$ be K\"ahler manifolds. When the manifolds are
irreducible Hermitian symmetric spaces of compact type equipped with
canonical K\"ahler-Einstein metrics, we always normalize the
canonical metrics to have the maximum holomorphic sectional
curvature $+2$ as that for the projective space equipped with the
Fubini-Study metric.
 Let
$(\lambda_1,\cdots,\lambda_v)$ and $(\mu_1,\cdots,\mu_m)$ be two
sets of positive real numbers. Let $U\subset M$ be a connected open
subset and $G_j: U\ra M_j'$ be a non-constant holomorphic map for
each $j$. We write
$\omega_{M,G,\lambda}:=\omega_M+\sum_{j=1}^{v}\lambda_jG_j^*\omega_{M'_j}$
for the modification of $\omega_M$ by $(G,\lambda):=(G_1,\cdots,
G_v;\lambda_1,\cdots,\lambda_v)$ over $U$. We are interested in the
 question: When are
there holomorphic maps $F_l: U\ra M_l$ 
for each $l \in \{1,\cdots, m\}$ such that
$\omega_{M,G,\lambda}=\sum_{l=1}^{m}\mu_l F_l^*\omega_{M_l}$? What
is the rigidity  phenomenon in this setting? There is a case where
we may not get anything interesting, due to the obvious cancelation.
In fact, as shown in Example 6.2, when we do not have the following
number theoretic property:
\begin{equation}\label {james-10}
\hbox{span}_{\mathbb Q_0^+}\{\lambda_j\}_{j=1}^v\cap
\hbox{span}_{\mathbb Q_0^+}\{\mu_l\}_{l=1}^{m} =\{0\},
\end{equation}
 we can easily construct
examples where we do not have any good rigidity and global extension
properties. Here we write ${\mathbb Q}_0^+$ for the set of
non-negative rational numbers.

A  main result of this paper is to provide the following rigidity
theorem, under the  the needed  number theoretic condition
(\ref{james-10}) for the conformal factors.

\begin{theorem}\label{kahler}
Let $h(z, \xi)=a_0+\sum_{|\a|,|\b|>0}a_{\a\b}z^{\a}\xi^{\b}$ be an
irreducible holomorphic polynomial over $\mathbb{C}^{2n}$ for $n
\geq 1$ with $h(z,\-{z})>0$.
Let $M$ be a simply connected  $n$-dimensional complex manifold
(not necessarily complete)
 with a real analytic K\"ahler metric $\omega_M$. Assume that
 there exists a holomorphic coordinate chart  ($U \subset X,\ \phi)$
 with $\phi(U)=V\ni 0$, a connected open subset in ${\mathbb C}^n$, such that
 $(\phi^{-1})^*\omega_M
 = \sqrt{-1}\partial \bar\partial \log h(z, \bar z)$. Let $M_l={\mathbb P}^{N_l}, M'_j={\mathbb P}^{N_j'}$
  be complex projective spaces 
equipped with the  Fubibi-Study   metrics $\omega_l, \omega'_j$,
respectively,  for $1 \leq l \leq m, 1 \leq j \leq v.$
 Suppose that $F_l: U\subset M \rightarrow {\mathbb P}^{N_l}, G_j: U\subset M \rightarrow {\mathbb P}^{N_j'}$ are non-constant holomorphic maps over $U$ for each $l,j$.
 Assume that
\begin{equation}\label{kahler isometry}
\omega_M = \sum_{l=1}^m \mu_l F_l^* \omega_l - \sum_{j=1}^v
\lambda_j G_j^* \omega'_j~~\text{over}~~U,
\end{equation}
for 
$\mu_l, \lambda_j \in \mathbb{R}^+$ which satisfy the property in
(\ref{james-10}). Then $F_l$ and $G_j$  extend to global holomorphic
immersions $\tilde{F_l}:\ M \ra {\mathbb P}^{N_l}$ and
$\tilde{G_j}:\ M\ra {\mathbb P}^{N_j'}$, respectively. Moreover
$\tilde{F_l}^*\omega_{l}=m_l\omega_M,
\tilde{G_j}^*\omega'_{j}=n_j\omega_M$ with $m_l, n_j\in {\mathbb N}$
satisfying  the identity:
\begin{equation}\label{james-0008}
1=\sum_{l=1}^m\mu_l m_l-\sum_{j=1}^v\lambda_j n_j.
 \end{equation}
\end{theorem}

In Theorem \ref{kahler}, when $M$ is not simply connected, then one
can still conclude that the maps extend along any path inside $M$
initiated from a point in $U$.  Theorem \ref{kahler} applies
immediately when $(M,\omega_M)$ is the projective space equipped
with the Fubini-Study metric. Indeed, the special case of Theorem
\ref{kahler} with $M$ being ${\mathbb P}^1$  can be  applied with
the minimal rational curve theory  on Hermitian symmetric spaces of
compact type to yield the following:

\bigskip
\begin{theorem}\label{maintheorem-local} Let  $(M,\omega_M), (M_l, \omega_{M_l}), (M'_j, \omega_{M_j'})$ be irreducible
compact Hermitian symmetric spaces
of compact type 
equipped with the canonical  K\"ahler-Einstein  metrics $\omega,
\omega_l, \omega'_j$ for $1 \leq l \leq m, 1 \leq j \leq v$,
respectively. Let $\mu_l, \lambda_j>0$ be a set of real numbers  for
$1\le l\le m$ and $1\le j\le v$, satisfying (\ref{james-10}).
Let $F_l: U\ra M_l$ and $G_j: U\ra M_j'$ be  non-constant
holomorphic maps for each $l, j$ such that
\begin{equation}\notag 
\omega_M=\sum_{l=1}^{m}\mu_lF^*_l\omega_{M_l}-\sum_{j=1}^{v}\lambda_j
G^*_j\omega_{M'_j}\ \ \hbox{over}\ \ U,
\end{equation}
where $U\subset M$ is a connected open subset.
Then for any $j\in \{1,\cdots,v\} \
\hbox{and}\ \ l\in \{1,\cdots, m\}$, $F_l$ and $G_j$  extend to
global holomorphic embeddings $\tilde{F_l}:\ M\ra M_l$ and
$\tilde{G_j}:\ M\ra M_j'$, respectively. Moreover
$\tilde{F_l}^*\omega_{M_l}=m_l\omega_M,
\tilde{G_j}^*\omega_{M'_j}=n_j\omega_M$ with $m_l, n_j\in {\mathbb
N}$ satisfying the equation (\ref{james-0008}).
\end{theorem}


Notice that in Theorem \ref{kahler}, the extended maps may not be
one to one, while in Theorem \ref{maintheorem-local}, they all are
embeddings. Also in these two theorems, when $v=0$, it is understood
that there are no mappings $G_j$. As one sees in Example 6.2, the
assumption (\ref{james-10}) on the conformal factors
 in Theorem \ref{kahler} and Theorem
\ref{maintheorem-local} is to avoid the cancelation that destroys
good rigidities, which together with (\ref{james-0008}) is necessary
and sufficient for Theorem \ref{kahler} and Theorem
\ref{maintheorem-local} to hold.
This phenomenon is similar  to the study of CR mappings:
Non-existence of negative conformal factors corresponds to the
strongly pseudoconvex case, while the investigation for the case of
 conformal factors with a mixed sign resembles to  the Levi non-degenerate situation
with positive signature (see [Hu2], [BH], for instance).


Theorem \ref{kahler} and Theorem \ref{maintheorem-local} follow
trivially from the Calabi theorem, when  all conformal factors
are positive and rational (or at least rational up to a common
factor). In this setting, one can  apply the Nakagawa-Takagi
theorem, the Veronese and Segre embeddings, to immediately reduce
them to the study a local holomorphic isometric embedding into a
projective space equipped with a multiple of the Fubini-Study
metric. Hence, the results follows immediately from the classical
Calabi theorem in [C].
(See Mok \cite{M4} and Yuan \cite{Y}.)
However,
by Calabi [p. 23, C] (see also Remark 6.1), when $\mu_1/\mu_2$ is
not a rational number, then the K\"ahler manifold $({\mathbb P}^n,
\mu_1 \omega_n)\times ({\mathbb P}^m, \mu_2 \omega_m)$ can not be
isometrically embedded into the Hilbert projective space equipped
with any multiple of the Fubini-Study metric $({\mathbb P}^{\infty},
\mu \omega_\infty)$. (Here, $\mu>0$ is a real number.)

In [C], Calabi also considered the existence problem for isometric
embeddings between complex manifolds with indefinite K\"ahler
metrics. Since in this setting, there is always a cancelation even
inside the metric part itself, it seems hard to get any rigidity and
global extension result.

A  first main   step in our proof of Theorems \ref{kahler}  is to
obtain the Nash-algebraicity for both $F_l$ and $G_j$ under the
assumption of (\ref{james-10}).
We then prove a single-valuedness for algebraic functions satisfying
a certain transcendental equation by employing the monodromy
argument and the Puiseux expansion for multi-valued algebraic
functions.
The last step is to use what is obtained in previous steps, as well
as, the geometry of minimal rational curves over the source
manifold, to derive 
Theorem
\ref{maintheorem-local}.


\medskip

The organization of the paper is as follows. In Section 2, after
collecting the basic properties of algebraic functions, we prove a
certain induction property and boundedness for algebraic functions
satisfying a certain transcendental functional identity. In Section
3 and Section 4, we prove the algebraicity and single-valuedness,
respectively, for germs of holomorphic functions satisfying a
certain transcendental functional identity. At the end of the section 4, we give the proof of Theorem \ref{kahler}. In Sections 5, we
apply the results derived in Section 3 and 4, to prove Theorem
\ref{maintheorem-local}. We end up our
paper with examples which show our number theoretic condition
(\ref{james-10}) on the conformal factors are more or less necessary
and sufficient conditions for our results to hold.

\medskip
{\bf Acknowledgement}: The authors would like to express their
gratitude  to  N. Mok and S.-C. Ng for many very valuable
discussions related to this paper. The  authors  would like also to
thank J.-M. Hwang for several helpful communications related to this
work in January, 2011.

\bigskip

\medskip

\section{ Algebraic functions}

In this section, we start by first recalling some basic properties
for algebraic functions. Then we will prove several lemmas to be
used in Section 3 and Section 4.
Our basic tools  are the monodromy argument and the Puiseux
expansion for multi-valued algebraic functions.


A holomorphic function $f$ over $U\subset {\CC}^n$ is called a
holomorphic Nash-algebraic, or simply holomorphic algebraic
function, if there is an irreducible polynomial $P(z,X)$  in $X$
with coefficients in polynomials of $z$ such that $P(z,f(z))\equiv
0$ over $U$. $f$ extends  to multi-valued a function, still denoted
by $f$, over ${\CC}^n$, some branches of which may blow up along a
complex analytic variety of codimension one. $P(z,X)$ is called a
minimum polynomial of $f$. All branches of $f$ are bounded over any
compact subset in ${\CC}^n$ if and only if the leading coefficient
of $P(z,X)$ can be made to be $1$. For our purpose in this article,
we are mainly concerned with algebraic functions whose minimum
polynomials have leading coefficient 1.

 Let $H$ be an
algebraic function in $U\subset \mathbb{C}^n$ with its minimum
polynomial of leading coefficient $1$: $P(z, Y) = Y^d+ a_1(z)
Y^{d-1} + \cdots + a_d(z)$, where $d \geq 1, a_1, \cdots, a_d$ are
polynomials in $z$ and $P(z, Y)$ is an irreducible polynomial in
$(Z,Y)$. Then there are polynomials $A(z, Y), B(z, Y)$ and $p(z)$
such that

$$A(z, Y)P(z; Y)+ B(z, Y)\frac{\partial P}{\partial Y}(z; Y) = p(z).$$
Let $E \subset \mathbb{C}^n$ be the affine algebraic hypersurface
defined by $p(z)=0$. Notice that any point $z_0$ with $p(z_0)\not
=0$ is  a  regular point for the algebraic function $H$, namely, a
point $z_0$  where for any $Y$ with $P(z_0,Y)=0$ it holds that
$\frac{\partial P}{\partial Y}(z_0, Y)\not =0.$
 Then any
 branch of $H(z)$ can be holomorphically continued along a curve $\gamma \subset \mathbb{C}^n \setminus E$.
  Let $E_0 \subset E$ be an irreducible component of $E$. We say that $E_0$ is a branching variety of $H$ if for a generic smooth point $p_0 \in E_0$
  and a sufficiently small ball $B_{p_0}$ centered at $p_0$, $\pi^{-1}(B_{p_0} \setminus E_0)$  has less than $d$
  pieces of connected components.
  Here, letting $W \subset \mathbb{C}^n \times \mathbb{C}$ be defined by $P(z, Y)=0$, then $\pi$ is a branched covering map
   from $W$ to $\mathbb{C}^n$ giving by the natural projection map. Since $p_0 \in E$ is a smooth point, after
   a holomorphic change of coordinates, we assume that $p_0=0$ and $E_0$ is defined by $z_n=0$.
   Being generic, we mean that $p_0=0$ is not contained in any other component of $E$.
   Let $\gamma(t)=(0, \cdots, 0, \epsilon e^{2\pi \sqrt{-1}t})$ for $0 \leq t \leq 1$, with $0 < \epsilon << 1$ sufficiently small.
   Then for a small ball $B_0$ centered at 0, any loop in $B_0 \setminus E_0$ is homotopic to $k \gamma$ with $k \in
   \mathbb{Z}$. Any simple loop homotopic to $\gamma$ is called a
   basic loop around $E_0$ near $p_0=0$.
    Now, write $p^*=(0, \cdots, 0, \epsilon)$. Then $E_0$ is a branching hypervariety for $H$ if and only if for some holomorphic branch
     $H_1$ of $H$ in a neighborhood of $p^*$, when we holomorphically continue $H_1$ along $\gamma$ one round,
     we will obtain another branch $H_2(\not= H_1)$ of $H$ near $p^*$. When $\frac{H_2}{H_1} = h_{12}$ is a constant,
     then we have some $d_0 \in \mathbb{N}$, such that, $h_{12}^{d_0}=1$.
     The smallest such a $d_0$ is called a period. In this case, we call $E_0$ a simple cyclic branching
     hypervariety for $H$ with respect to the
      branch $H_{1}$.
     Apparently, the simple cyclicity of $H$ along $E_0$ is independent of the choice of the
     generic point $p_0$.
  Furthermore, if there does not exist a polynomial function $h$ and a natural number $n_1$ such that
     $H=h^{\frac{1}{n_1}}$,
      then there must be an $n_2 \in \mathbb{N}$ such that $H^{n_2}$ has branching varieties, none of which is a simple cyclic branching
      variety with respect to any
       branch of  $H^{n_2}$.
       We also recall the Puiseux expansion for  branches obtained by continuing $H_1$ near $p_0=0$,
       which will be the basic tool for us to deal with
       multi-valued holomorphic functions:

$$H_1=\sum_{i=0}^\infty a_i(z') z_n^{\frac{i}{N_0}}$$ with $N_0(\geq 2) \in \mathbb{Z}$ and $a_i(z')$ being
holomorphic  near $0'$. Here $z'=(z_1, \cdots, z_{n-1})$.




\medskip

\medskip

The following lemma is essential in the later induction argument.

\begin{lemma}\label{branch} Let $U \subset \mathbb{C}^{n}$ be a connected open subset.
Let $\vec{h}(z)=(h_1, \cdots, h_k) (z) \in \mathbb{C}^k$ be a row of irreducible distinct polynomials. 
Let $\vec{r} \in ({\RR}\setminus \{0\} )^k$. Let $H_1(z), \cdots,
H_k(z)$ be non-constant holomorphic algebraic functions defined over
$U$
 such that for  a certain branch of the power function
 $(H_\a)^{\mu_\alpha}$ for each $\a$, and for a certain constant $c_0$, we have
\begin{equation}\label{equa}
\vec{h}(z)^{\vec{r}} = c_0\prod_{\alpha=1}^k H_\alpha(z)^{\mu_\alpha} ~~\text{on}
~~ U',
\end{equation}
for 
$\mu_1, \cdots, \mu_k \in {\mathbb R}\sm \{0\}.$
\ Here $U'\subset U$ is a subdomain such that
$\vec{h}(z)^{\vec{r}}=h_1^{r_1}\cdots h_k^{r_k}$ has a well-defined
holomorphic branch.
 Assume all the branches of $H_\alpha$'s can only have zero or
 points of indeterminancy
in the variety defined by the union of of the zeros of $h_j$ for
each $j$.
%
 If $H_1$ has an irreducible
non-cyclic branching hypervariety with respect to a certain branch,
then there exists $n^+_1, n_2, \cdots, n_k \in \mathbb{Z}$ with
$n^+_1>0$, such that
\begin{equation}\label{com}
\mu_1 n^+_1 = \sum_{\alpha=2}^k \mu_\alpha n_\alpha.
\end{equation}
Furthermore,
we have
\begin{equation}\label{equa2}
\vec{h}(z)^{n^+_1 \vec{r}} =c_0^{n^+_1} \prod_{\alpha=2}^k \left(
H_1^{n_\alpha} H_\alpha^{n^+_1}\right)^{\mu_\alpha} .
\end{equation}
\end{lemma}

In the above lemma and the rest of the paper, for an algebraic
function $\phi$, we define the zero set  (points of indeterminancy,
respectively) of $\phi$ to be that defined by the zeroth order  term
(the coefficient of the highest order term, respectively) in a
minimal polynomial of $\phi$.

\medskip
 {\it Proof of Lemma \ref{branch}:} Assume without
loss of generality that $0$ is in $ U$  and is a regular point for
$H_\a$ for each $\a$.
Since $H_1$ has an irreducible non-cyclic branching hypervariety $E$
with respect to a certain branch of $H_1$. Let $p^*(\not\in E)$ be
very close to a generic smooth point $p^*_0$ of $E$. In what
follows, we assume that $E$ is defined by the zero of $h_i(z)$ for a
certain $i$ or a small neighborhood of $p^*_0$ does not cut any zero
of $h_\a's$.
Let $\gamma$ be a basic loop around $E$ near $p^*_0$ with
$\gamma(0)=p^*$.
Assume that $\gamma$ stays in a sufficiently small ball centered at
$p^*_0$. Assume $(H_1)_1$ is a holomorphic branch of $H_1$ at $p^*$
and when we continue holomorphically $(H_1)_1$ along $\gamma$ one
round, we get a new holomorphic branch $(H_{1})_2$ near $p^*$ with
$\chi_{11}=\frac{(H_{1})_2}{(H_1)_1} \not= \text{constant}$. Choose
$\gamma_0$, a simple curve connecting 0 to $p^*$ (that does not cut
the zero of any $h_j$  and the union of branching varieties of $H_l$
for each $l$) such that we get $(H_{1})_1$ near $p^*$ when we
continue $H_1$ along $\gamma_0$ from its original holomorphic value
near $0$ . 
We still have by the uniqueness of real analytic functions:
\begin{equation}\label{branch1-lemma}
\vec{h}^{\vec{r}} = c_0(H_1)_1 ^{\mu_1} 
\cdots (H_k)_1^{\mu_k} ~~~\text{near}~~~ p^*.
\end{equation}
Here $(H_\alpha)_1^{\mu_\alpha}$ 
 is a certain branch of the multi-value
functions $e^{\mu_\alpha \log (H_\alpha)_1}$ 
 near $p^*$ for any $\alpha$ with $1 \leq \alpha \leq k$.
 Now we holomorphically continue
(\ref{branch1-lemma}) along $\gamma$ to arrive at

\begin{equation}\label{branch2-lemma}
\vec{h}^{\vec{r}} =c_0'(H_1)_2^{\mu_1}
\cdots (H_k)_2^{\mu_k} ~~~\text{near}~~~ p^*,
\end{equation}
since  $h^{r_i}_i$ is at most cyclically branching at $p$ for each
$1\leq i \leq k$.
From (\ref{branch1-lemma}) and (\ref{branch2-lemma}), we get
\begin{equation}\notag
 c_{1}= \chi_{11}^{\mu_1} 
\cdots \chi_{1k}^{\mu_k}~~~\text{near}~~~ p^*
\end{equation}
with $c_{1}(\not= 0)$ a certain constant, reflecting how we choose
the branches of the multi-valued functions. Hence
$$\chi_{11} = c'_{1} \chi_{12}^{-\mu_2 / \mu_1}\cdots \chi_{1k}^{ -\mu_k / \mu_1}\ \hbox{near  } p^*.$$ Here, $c'_{1}(\not= 0)$ is a certain constant.
Since $\chi_{11}$ is not a constant, we can find a complex line $L$
with a linear coordinate $\xi$ such that $\tau = \chi_{11}(\xi)
\not= \text{constant}$. Here and in what follows, we will write
$\chi_{11}(\xi)$ for $\chi_{11} \big|_L$. Hence

$$\tau = c''_{1} (\chi_{12} \circ \chi_{11}^{-1}(\tau))^{-\mu_2 / \mu_1}\cdots (\chi_{1k} \circ \chi_{11}^{-1}(\tau))^{-\mu_k / \mu_1}.$$
Now, after a holomorphic continuation, we can assume that the above
holds for certain branches near the origin.
 Since $\chi_{1\alpha} \circ \chi_{11}^{-1}(\tau)$ is algebraic for $2 \leq \alpha \leq k$, we have the Puiseux expansion
 near
 0:
$$\chi_{1\alpha} \circ \chi_{11}^{-1}(\tau) = \sum_{i \geq i_\alpha} a_{\alpha i} \tau^{i / N_0},\ a_{\alpha i_\alpha}\not= 0.$$
Thus, we get (\ref{com}) by comparing the exponent of the lowest
degree term in the expansion in $\tau$. 
%
Taking $n^+_{1}$-th power in (\ref{equa}) and applying (\ref{com}),
we get (\ref{equa2}).
 $\endpf$
%
\medskip





In the following, for an open subset $U\in {\CC}^n$, we write
$\hbox{conj}(U):=\{z,\-z\in U\}$. The following is a key lemma for
our proof of Theorem \ref{singleness123}.

\begin{lemma}\label{leading coeff} Let $U_1$ be a   connected open
subset in ${\CC}^n$, and let
$\phi(z,\-{z})$ be real analytic in $(z,\-{z})$ over ${\CC}^n$ with
$\phi(z,\-{z})>0$ for any $z\in {\CC}^n$. Assume that $\phi(z,\xi)$
is holomorphic and algebraic in $U_1\times \hbox{conj}(U_1)$.
 Suppose that
 $F_l(z): U_1 \subset \mathbb{C}^{n_1} \rightarrow
\mathbb{C}^{N_l}, G_j(z): U_1 \subset \mathbb{C}^{n_1} \rightarrow
\mathbb{C}^{N'_j}$ non-constant holomorphic algebraic maps. Write
$\-{F}(\xi):=\-{F(\-{\xi})}$ for $\xi\in \hbox{conj}(U_1)$.
Suppose they satisfy the following
identity

\begin{equation}\label {new-01}
\prod_l (1+ F_l(z) \cdot \overline{F_l}(\xi))^{\mu_l} = \phi(z, \xi)
\cdot \prod_j (1+ G_j(z) \cdot \overline{G_j}(\xi))^{\nu_j}
~~\text{over} ~~
 U_1 \times \hbox{conj}(U_1),
 \end{equation}
 where $\mu_j,\nu_j$ are positive real numbers satisfying
 (\ref{james-10}).
%
Then the following holds:

\begin{enumerate}
\item Each component of $F_l(z)$ and $G_j(z)$ has a minimal polynomial with
leading coefficient 1. In particular, for any compact subset $K_1 \subset \mathbb{C}^{n_1}$,
there is a constant $C_{K_1}$ such that $|h^*|\le C_{K_1}$ holds over $K_1$,
for any branch $h^*$ of a function component $h$ from $G_j$ or
$F_l$.

\item Any branch of
$1+ F_l(z) \cdot \overline{F_l}(\xi)$ and $1+ G_j(z) \cdot
\overline{G_j}(\xi)$ can only have zeros on the hypervariety
consisting of the  zeros and points of inderterminancy of  $\phi(z,
\xi)$. (When $\phi(z,\xi)$ is not a polynomial, as mention before,
the zero set ( the set of points of indeterminancy) of $\phi(z,\xi)$
is the set defined by the zeroth order term (the coefficient of the
highest order term, respectively) in a minimal polynomial of
$\phi(z,\xi)$)
\end{enumerate}
\end{lemma}

{\it Proof of Lemma \ref{leading coeff}}: 
For simplicity of notation, assume that the first component of
$F_1$, denoted by $h_0(z)=F_{1,1}(z)$, has a minimial polynomial
with leading coefficient $a(z)\not= \hbox{constant}$. Let $a_0(z)$
be an irreducible factor of $a(z)$, and let $K_0(z)=a(z)h_0(z)$.
Then $K_{0}$ is an algebraic function with a minimum polynomial of
leading coefficient $1$. Write $\{h_1,\cdots, h_k\}$ for the other
function components in $F$ and $G$ whose respective
leading coefficients $a_1,\cdots,a_k$   have the prime factor $a_0(z)$. 
Define similarly $K_{\alpha}(z)=a_\alpha(z)h_\alpha(z)$ for
$\alpha=1,\cdots,k$. Choose a point $p_0$ such that (i) $p_0$ is a
smooth point of the zero set of $a_0(z)$,
(ii) the leading coefficients of the rest components
 are not zero at $p_0$ and the other components of $a_1, \cdots, z_k$ are not zero near $p_0$, and (iii) the zero set of $a_0$ near $p_0$, which may or may not be a branching variety itself,
   does not cut any other branching variety of  components of $F_l,
 G_j$.
  Assume $p_0=0$ to simplify the notation.
After a local change of
 coordinates, we may assume that $a_0(z)=z_n$. Now, by the analytic extension we can assume (\ref{new-01}) with $\xi=\-{z}$ holds for $(z,\-{z})$
 near $p_0=0$
 for certain branches of $F_l$ and $G_j$. For $F_{1,1}$, the branch is chosen so that it is unbounded near
 $p_0$. (From the relation formula between
roots and coefficients, and using the irreducibility of minimal
polynomials, we can see that not all branches of $F_{1,1}$ can be
bounded near $p_0$.  Hence, we can pick the unbounded one to fit our
consideration here.)
 Consider the Puiseux
 expansion of the corresponding branch $K_\alpha$ near $p_0$ for $0 \leq \alpha \leq k$, which we still write as $K_\alpha$ to simplify
 the notation: 

 $$K_\alpha=
 \sum_{i=i_{\alpha} \geq 0}^{\infty}b_{\alpha, i}(z')z_n^{\frac{i}{n_\alpha}}$$
 with $b_{\alpha, i_\alpha}(z')\not \equiv 0$.
 Then
 $$|h_\alpha|^2=\big| z_n \big|^{2(\frac{i_\alpha}{n_\alpha} -k_\alpha)} |b_{\alpha, i_\alpha}(z')|^2+
 o(| z_n |^{2(\frac{i_\alpha}{n_\alpha}
 -k_\alpha)} ),$$
 with $k_\alpha$ a certain natural number.

We  easily see that there exist rational numbers $r_l \leq 0, r'_j
\geq 0$ such that
$$|z_n|^{(\sum_{l}r_l\mu_l+\sum_{j}r_j'\lambda_j)}R(z,\-{z})= \phi(z, \bar z)$$
near $p_0=0$. Notice that $r_1\le 2\left(\frac{i_1}{n_1}
-k_1\right)<0$, for, as we mentioned,  the branch $h_0$ of $F_{1,1}$
is chosen to be unbounded near $p_0$. (Similarly, we have
$r_\alpha\le 0$ when the associate component is from $H$'s and
$r'_\alpha\ge 0$ when the associate component  is from $G$'s.)
 Here for a generic fixed choice of $z'$ near $0$, when
$z_n\ra 0$, $|R(z,\-{z})|$ and $|R(z,\-{z})|^{-1}$ remain to be
bounded. Hence, we have $\sum_{l}r_l\mu_l+\sum_{j}r_j'\lambda_j=0$,
contradicting the hypothesis in (\ref {james-10}). Thus Lemma
\ref{leading coeff} (1) is proved.

\medskip

To prove the second part of the lemma, 
suppose that $p_0$ is not a zero nor a point of indeterminancy of
$\phi(z, \xi)$, but is a zero for a certain branch $K$ of $H_l:=1+
F_l(z) \cdot \overline{F_l}(\xi)$ or $Q_j:=1+ G_j(z) \cdot
\overline{G_j}(\xi)$. Consider the variety $Z_{p_0}$ defined by the
coefficient polynomial of the degree zero term in the minimum
polynomial of $K$. Write $X$ for the union of the branching
varieties of $H_l,\ Q_j$. If $Z_{p_0}\not \subset X$, we move $p_0$
slightly such that $p_0\not\in X$. Otherwise, we slightly move $p_0$
such that $p_0$ is a smooth point of $X$. Then we consider the
Puiseux expansion of the involved branches of $G_j, H_l$ near $p_0$
and compare the vanishing orders near $p_0$. As in Lemma
\ref{leading coeff} (1), we finish the proof by applying
(\ref{james-10}). $\endpf$

\medskip

We next formulate  two elementary lemmas, which will be used in our
induction argument later.

\begin{lemma}\label{combo}
Let $(a_1, \cdots, a_p)$ and $(b_1, \cdots, b_p)$ be  $p$-tuples
consisting of positive numbers with $p\ge 1$. Then there exists
$i_0$ such that

$$\frac{a_i}{a_{i_0}} \geq \frac{b_i}{b_{i_0}}~~\text{for~all}~~i.$$
\end{lemma}

{\it Proof of Lemma \ref{combo}}: $i_0$ can be simply  chosen so
that  $\frac{a_i}{b_i}$ achieves the minimum value when $i=i_0$. $\endpf$

\medskip
\begin{lemma}\label{new-lemma} Let $F$ and $G$ be two vector valued
holomorphic factions near $0\in {\CC}^n$. For any two non-negative
integers $m_1, m_1$, there is a vector valued holomorphic function
$H$ near $0\in {\CC}^n$ such that $(1+|F|^2)^{m_1}\cdot
(1+|G|^2)^{m_2}=1+|H|^2.$
\end{lemma}

{\it Proof of Lemma \ref{new-lemma}}: When $m_1,m_2=0$, simply take
$H=0$. The other case  follows from the following elementary algebra
formula: Let $a=(a_1,\cdots,a_k)$ and $b=(b_1,\cdots,b_k)$. Then
$$(1+|a|^2)(1+|b|^2)=1+|a|^2+|b|^2+\sum_{j,l=1}^k |a_jb_l|^2. ~\endpf$$


\section{Algebraicity for holomorphic maps satisfying  transcendental functional equations}
In this section, we prove the algebraicity for germs of holomorphic
functions satisfying a certain transcendental relation. 

\begin{proposition}\label{gen alg}
Let $F_l: B\subset \mathbb{C}^n \rightarrow \mathbb{C}^{N_l}, G_i:
B\subset \mathbb{C}^n \rightarrow \mathbb{C}^{N'_i}$ be holomorphic
maps defined over a small ball $B$ centered at the origin for $1
\leq l \leq m, 1\leq i \leq v$ with $F_l(0), G_j(0)=0$. Let
$\phi(z, \xi): B\times B \rightarrow \mathbb{C}\sm \{0\}$ be a
holomorphic algebraic function. Suppose that
there exists $\mu_l, \lambda_i \in \mathbb{R}^+$ such that
\begin{equation}\label{log}
 \phi(z, \xi) \cdot \prod_{i=1}^v \left(1+
G_i(z) \cdot \bar G_l(\xi)\right)^{\lambda_i} = \prod_{l=1}^m
\left(1+ F_l(z) \cdot \bar F_l(\xi)\right)^{\mu_l} ~~\text{over}~~B
\times B.
\end{equation}
Then there exist  algebraic functions $\wh{F}_l(z, X_1, \cdots,
X_v)$ holomorphic for  $(z, X_1, \cdots, X_v)\approx 0$ for all $1
\leq l \leq m$ such that $F_l(z)=\wh{F}_l(z, G_1(z), \cdots,
G_v(z))$ for $z\approx 0$. In particular, when $v=0$, then each
$F_l$ is algebraic. (Notice here, there is no restriction on $\mu_j$
except that they are positive.)
\end{proposition}

{\it Proof of Proposition \ref{gen alg}}: Write
$D_\delta={\partial\over\partial z_\delta}$ and for
$\alpha=(\alpha_1,\cdots,\alpha_n)\in ({\mathbb Z_0^+})^n$, write
$D^{\alpha}={\partial^{|\alpha|}\over\partial z_1^{\alpha_1}\cdots
\partial z_{n}^{\alpha_n}}$.
Here, we  write ${\mathbb Z}_0^+$ for the set of non-negative
integers. Applying   the
logarithmic differentiation to (\ref{log}), we get for any $1 \leq
\delta \leq n$ and $z,\xi\approx 0$ the following:
\begin{equation}\label{james-007-2}
\sum_{l=1}^{m}{\mu_l D_\delta(F_l)(z)\cdot \-{F_l}(\xi)\over
1+F_l(z)\cdot \-{F_l}(\xi)}=\wh{ \phi}_\delta(z,\xi) +
\sum_{i=1}^{v}{\lambda_i D_\delta(G_i)(z)\cdot \-{G_i}(\xi)\over
1+G_i(z)\cdot \-{G_i}(\xi)},
\end{equation}
where $ \wh{\phi}_\delta(z,\xi)=D_\delta \log
\phi (z, \xi)$ is Nash algebraic in $(z,\xi)$ by
the assumption.
Write
$$\chi=(\chi_1,\cdots,\chi_N)=(\sqrt{\mu_1}F_1,\cdots,\sqrt{\mu_m}F_m)$$
with $N=N_1+\cdots N_m.$ Then, we can rewrite (\ref{james-007-2}) as
follows:
\begin{equation}\label{james-007-3}
\sum_{l=1}^{N} D_\delta(\chi_l)(z)\cdot \-{\chi}_l(\xi) +
H_\delta(z,\-{\chi}(\xi))= \Phi_\delta(z,\xi, \bar G_1(\xi), \cdots, \bar G_v(\xi)),
\end{equation}
where $ \Phi_\delta(z,\xi, X_1(\xi), \cdots, X_v(\xi))$ is a Nash
algebraic function in $\xi$ and $X_i(\xi)$ for all $i$. Now,
differentiating (\ref{james-007-3}), we get for any $\alpha$ the
following

\begin{equation}
\label{james-007-4} \sum_{l=1}^{N} D^{\alpha}(\chi_l)(z)\cdot
\-{\chi_l}(\xi) +
H_{\alpha}(z,\-{\chi}(\xi))=\Phi_{\alpha}(z,\xi, \bar G_1(\xi), \cdots, \bar G_v(\xi)).
\end{equation}
Here for any $\alpha$ and for any fixed $z$,
$H_{\alpha}(z,\-{\chi}(\xi))$ is rational in $\-{\chi}$, and has no
constant  and linear terms in the Taypor expansion with respect to
$\-{\chi}$. Also, $\Phi_{\alpha}(z,\xi, X_1, \cdots, X_v)$ 
is Nash algebraic in
$\xi, X_1, \cdots, X_v$ for any fixed $z$. 

We emphasize that, the advantage in the present situation is that
$z, \xi$ are totally independent variables. Now, let ${\mathcal
L}:=\hbox{Span}_{\mathbb C}\{ D^{\alpha}(\chi(z))|_{z=0}\}_{|\a|\ge
1}$ be a vector subspace of ${\mathbb C}^{N}$. Let
$\{D^{\alpha_i}(\chi(z))|_{z=0}\}_{i=1}^{\tau}$ be a basis for
$\mathcal L$. Then for a small open ball $B_0$ centered at $0$ in
${\mathbb C}^n$, $\chi(B_0)\subset {\mathcal L}.$ Indeed, for any
$z\approx 0$, we have from the Taylor expansion that
$\chi(z)=\chi(0)+\sum_{|\alpha|\ge 1}{D^{\alpha}(\chi)(0)\over\alpha
!}z^{\alpha}=\sum_{|\alpha|\ge 1}{D^{\alpha}(\chi)(0)\over\alpha
!}z^{\alpha}\in {\mathcal L}.$

Now, let $\nu_i
~(i=1 \cdots,N-\tau)$ be a basis of the Euclidean orthogonal
complement of ${\mathcal L}.$  Then, we have

\begin{equation}
\label{james-007-06} \nu_{i} \cdot \-{\chi}(\xi)=0, \ \ i=1,\cdots,
N-\tau.
\end{equation}
Consider the system consisting of (\ref{james-007-4}) at $z=0$ (with
$\alpha=\alpha_1,\cdots,\alpha_\tau$) and (\ref{james-007-06}).
Since the Jacobian matrix of the functions in the left hand side of
the system with respect to $\-{\chi}$ at $0$ is
\begin{equation}
\begin{bmatrix}
D^{\alpha_1}(\chi(z))|_{z=0}\\
\vdots\\
\nu_{N-\tau}
\end{bmatrix}
\end{equation}
and is obviously invertible. Note that the left hand side of the
system of equations consisting of (\ref{james-007-4}) at $z=0$ (with
$\alpha=\alpha_1,\cdots,\alpha_\tau$) and (\ref{james-007-06}) is
Nash algebraic in $\xi$ and the right hand side is Nash algebraic in
$\xi, \bar G_1(\xi), \cdots, \bar G_v(\xi)$. By the algebraic
version of the inverse function theorem (see, for instance [Hu1]),
there exists Nash algebraic functions $\wh{F}_l(\xi, X_1, \cdots,
X_v)$ in all variables $\xi, X_1, \cdots, X_v$ for all $1 \leq l
\leq m$ such that $F_l(\xi)=\frac{\chi_l}{\sqrt{\mu_l}}(\xi)=
\wh{F}_l(\xi, G_1(\xi), \cdots, G_v(\xi))$ near $\xi=0$. This proves
Proposition \ref{gen alg}. $\endpf$

\bigskip
\begin{remark}
If one only has the following kind of
 equation
$$  \rho^*(F,\overline{F(z)})=k(z,\-{z})\rho(z,\-{z})$$ for a certain smooth function
$k(z,\-{z})$ with $\rho,\rho^*$ algebraic in their variables,
 then the differentiation
 can only be taken by CR vector fields  tangent to and along the manifold defined by $\rho(z,\-{z})=0$.
  The argument
 will be
 much more involved, and the algebraicity might  be only achieved under certain non-degeneracy assumptions.
 This is the situation that one encounters in
 studying CR maps. (See [Hu1] for instance).
\end{remark}

Let $F_l, G_j$ be holomorphic maps as in Proposition \ref{gen alg}
satisfying the transcendental equation (\ref{log}). It follows
 from Proposition \ref{gen alg} that there exists  algebraic functions $\wh{F}_l(z, X_1, \cdots, X_v)$ such that $F_l(z)= \wh{F}_l(z, G_1(z),
  \cdots, G_v(z))$ as germ of holomorphic functions at $z=0$ for all $l$. Let $N'=N'_1 + \cdots N'_v$ and let $g_1(z), \cdots, g_{N'}(z)$ be all
  components in $G(z)=(G_1(z), \cdots, G_v(z))$. Let $\mathfrak{R}$ be the field of  rational functions in $z$   and consider the field
  extension $$\mathfrak{F}=\mathfrak{R}( g_1(z), \cdots, g_{N'}(z) ).$$ Let $K$ be the transcendence degree of the field extension
   $\mathfrak{F} / \mathfrak{R}$. If $K = 0$, then each element in $\{g_1(z), \cdots, g_{N'}(z) \}$ is a Nash algebraic function. Hence each $F_l(z)$ is
   also algebraic for all $l$. Otherwise, by re-ordering the lower index, let $\mathcal{G}= \{g_1, \cdots, g_K\}$ be the maximal algebraic
   independent subset in $g_1, \cdots, g_{N'}$, and it follows that the transcendence degree of $\mathfrak{F} / \mathfrak{R}(\mathcal{G})$ is 0.
    In fact, for any $l>K$, there exists a minimal polynomial $P_l(z, X_1, \cdots, X_K, X)$ such that $P_l(z, g_1(z), \cdots, g_K(z), g_l(z))\equiv 0$
     and moreover, $\frac{\partial P_l(z, X_1, \cdots, X_K, X)}{\partial X}(z, g_1(z), \cdots, g_K(z), g_l(z)) \not\equiv 0$ in $U$, a small neighborhood of $0$,
     for otherwise, $P_l$ can not be a minimal polynomial of $g_l$. Now the vanishing of the partial derivatives for all $l$ forms a proper local
      complex analytic
     variety near $0$. Let $\gamma : [0,1]\ra U$ be a smooth simple curve with $\gamma(0)=0$ and $\gamma((0,1])$ does not curve the just mention variety.
     Applying the algebraic version of the existence and uniqueness part of the
     implicit function theorem,
      there exist a small connected open subset $U_0 \subset U$ with $0\in \-{U_0}$
      and  a holomorphic  algebraic function $\wh{g}_l$ in the neighborhood  $\wh{U_0}$ of
      $\{(z, g_1(z), \cdots, g_K(z)): z \in U_0\}$
      in ${\CC}^{n}\times {\CC}^{K}$,
       such that $g_l(z) = \wh{g}_l (z, g_1(z), \cdots, g_K(z))$ for any $z \in
       U_0$. (We can assume that $U_0$ is the projection of
       $\wh{U_0}$.)
        Substitute into $\wh{F}_l(z, G_1(z), \cdots, G_v(z))$,
       and still denote it, for simplicity of notation,  by $\wh{F}_l(z, g_1(z), \cdots, g_K(z))$ with
      $$\wh{F}_l(z, g_1(z), \cdots, g_K(z)) = \wh{F}_l(z, G_1(z), \cdots, G_v(z))\  \hbox{for}\  z \in
      U_0.$$
      Write $\wh{G}_j(z, g_1(z), \cdots, g_K(z) )$ for $\wh{g}$ when
      $j>K$ and $\wh{G_j}(z, g_1(z), \cdots, g_K(z) ) \ =g_j$ for $j\le
      K$.
Now, let $X=(X_1, \cdots, X_K)$ and replace 
$g_i(\xi)$ by $X_i$ for $1 \leq i \leq K$ in $F_l(\xi) =
\wh{F}_l(\xi, g_1(\xi), \cdots, g_K(\xi))$ and $G_j(\xi) =
\wh{G}_j(\xi, g_1(\xi), \cdots, g_K(\xi))$, in the following
quantity:

\begin{equation}
\sum_{l=1}^{m}{\mu_l \overline{D_{z_\delta}(F_l)(z)}\cdot {F_l}(\xi)\over
1+\overline{F_l(z)}\cdot F_l(\xi)}  - \sum_{j=1}^{v}{\lambda_i \overline{D_{z_\delta}(G_j)(z)}\cdot G_j(\xi)\over
1+\overline{G_j(z)}\cdot G_j(\xi)} - \frac{D_{\bar z_\delta} \overline{\phi}(\bar z,\xi)}{\overline{\phi}(\bar z, \xi)},
\end{equation}

Denote the new quantity by $\Phi_\delta (\bar z, \xi, X)$. We have
the following:

\begin{lemma}\label{albert}
$\Phi_\delta (\bar z, \xi, X) \equiv 0$ for any $1 \leq \delta\leq
n$, for any $z$ near 0 and any $(\xi, X) \in \wh{U_0}$.
\end{lemma}

{\it Proof of Lemma \ref{albert}:} Suppose not. Notice that
$\Phi_\delta (\bar z, \xi, X)$ is Nash algebraic in $(\xi, X)$ by
Proposition \ref{gen alg}. For a generic fixed $z=z_0$ near 0, since
$\Phi_\delta (\bar z_0, \xi, X) \not \equiv 0$, there exist
polynomials $A_{\alpha}(\xi, X)$ for $0 \leq \alpha \leq m$ with
$A_{0}(\xi, X) \not\equiv 0$ such that

$$\sum_{0 \leq \alpha \leq m} A_{\alpha}(\xi, X) \Phi_\delta^\alpha(\bar z_0, \xi, X) \equiv 0.$$
 As $\Phi_\delta (\overline{z_0}, \xi, g_1(\xi), \cdots, g_K(\xi))
\equiv 0$ for $\xi \in U_0$, then it follows that $A_0 (\xi,
g_1(\xi), \cdots, g_K(\xi)) \equiv 0$ for $\xi \in U_0$. This is a
contradiction to the assumption that $\{ g_1(\xi), \cdots, g_K(\xi) \}$ is
an algebraic independent set. $\endpf$

\medskip
We now are in a position to prove the following algebraicity result:

\begin{theorem}\label{algebraicity-1}
Let $(\lambda_1,\cdots,\lambda_v)$ and $(\mu_1,\cdots,\mu_m)$ be two
sets of positive real numbers satisfying the following condition:

$$ \hbox{span}_{\mathbb
Q_0^+}\{\lambda_j\}_{j=1}^v\cap \hbox{span}_{\mathbb
Q_0^+}\{\mu_l\}_{l=1}^{m}=\{0\}.$$
 Let $B$ be a small ball centered at the origin of ${\CC}^n$.
Let $F_l: B \rightarrow \mathbb{C}^{N_l}, G_j: B\ra
\mathbb{C}^{N'_j}$ be non-constant holomorphic maps with $F_l(0),
G_j(0)=0$ satisfying the transcendental equation:

\begin{equation}\label{new-02}
\prod_{l=1}^m \left(1+ F_l(z) \cdot \bar F_l(\xi)\right)^{\mu_l}
\cdot \prod_{i=1}^v \left(1+ G_i(z) \cdot \bar
G_l(\xi)\right)^{-\lambda_i} = \phi(z, \xi).
\end{equation}
  Suppose that $\phi(z,
\xi)$ is holomorphic and algebraic over $B\times B$.
Then $F_l$ and $G_j$ are Nash algebraic.
\end{theorem}

{\it Proof of Theorem  \ref{algebraicity-1}:} Keep the same notation
we set up above. After a unitary transformation  to simplify $F_l,
G_j$, if needed, we may assume in what follows that the component
functions in $F_l$ or $G_j$ are linearly independent over $\CC$.
Write
$$\Psi(z, \xi, X)=$$
$$ = \log \left(\prod_{l=1}^m \left(1+ \overline{F_l(z)} \cdot
\wh{F}_l(\xi, X)\right)^{\mu_l} \right) - \log \left(\prod_{j=1}^v
\left(1+ \overline{G_j(z)} \cdot \wh{G}_j(\xi, X)\right)^{\lambda_j}
\right) - \log \-\phi(\bar z , \xi).$$
Lemma \ref{albert} shows  that for any $z$ near 0, and $(\xi, X) \in
\wh{U_0}$ as defined there,
%
$\frac{\partial}{\partial \bar z_\delta} \Psi(z, \xi, X) =0. $
%
Notice that $\frac{\partial}{\partial z_\delta} \Psi(z, \xi, X) =0$
 as $\Psi ( z, \xi, X)$ is anti-holomorphic in $z$. Since $\Psi(0, \xi, X) =0$, then $\Psi(z, \xi, X) \equiv 0$ for $z\in B$ and  $(\xi, X) \in \wh{U_0}$.
Hence, we arrive at  the following identity,

\begin{equation}\label{albert111}
\prod_{l=1}^m \left(1+ \overline{F_l(z)} \cdot \wh{F}_l(\xi,
X)\right)^{\mu_l} = \overline\phi(\bar z , \xi) \prod_{j=1}^v \left(1+
\overline{G_j(z)}
 \cdot \wh{G}_j(\xi, X)\right)^{\lambda_j}\ \ (z,\xi, X)\in B\times \wh{U_0}.
\end{equation}

Here $\wh{F}_l$ and $\wh{G}_j$ are algebraic in their variables.
Notice that when $\wh{F}_l$ and $\wh{G}_j$ are independent of $X$,
then $F_l$ and $G_j$ are already algebraic. Hence, we will assume,
in the course of the proof of the theorem, that  one of the maps
from $\{\wh{F}_l\}$ or $\{\wh{G}_j\}$ depends on $X$. Let $h(\xi,X)$
be one of these. Let $E$ be the zero defined by the prime factors
involving only $\xi$ in the non-zero polynomial coefficients of a
minimal polynomial of $h(\xi,X)$. Then for a fixed $\xi\not \in B\sm
E$, $h(\xi,X)$ is not constant in $X$. Now, choose $X_\xi$ such that
$(\xi,X_\xi)\in \wh{U_0}$. Then for any unit vector $v$ in
${\CC}^K$, $h(\xi,t):=h(\xi,X_\xi+tv)$ is algebraic and holomorphic
in $t(\approx 0)\in \CC$. For a generic choice of the unit vector
$v$, $h(\xi,t)$ is a non-constant algebraic function in $t$.  Also
fix such a vector $v$. Now holomorphically continuing $h(\xi,t)$
along loops in $\CC$ avoiding its branch points ${\cal B}$, we get
multiple valued functions: $\{h_1(\xi,t),\cdots, h_m(\xi,t)\}$. If
for any compact subset $L$, $|h_j(\xi,t)|\le C_L$ for $t\in L\sm
{\cal B}$ with $C_L$ a constant depending only on $L$, then any
symmetric function of them is holomorphic over $L$. Certainly not
all symmetric functions of them can be uniformly bounded in ${\CC}$,
for, otherwise, by the Louville theorem, all symmetric functions of
$\{h_j(\xi,t)\}$ are constant and thus each $h_j$  is constant in
$t$, which contradicts to the assumption. Hence, in this setting, we
must have an unbounded branch near $t=\infty$. Let $t_0\in \CC$ be
such that some branches of $h$ are unbounded near $t_0$. For
simplicity of notation, assume $t_0=\infty$. Let $\gamma: [0, 1]
\rightarrow \mathbb{C}$ be a curve connecting $0$ to a point close
to $t_0=\infty$ and continue $h$ along $\gamma $ to get a branch,
which is  still denoted by $h$, near $\gamma(1)$. By what we argued
above, we can assume  that the continuation of such an $h$ along
curves near $\infty$ leads to unbounded multi-valued functions near
$t=\infty$.

For simplicity of notation, assume that the $h(\xi, X_\xi+tv)$
mentioned above is $\wh{F}_1(\xi, X_\xi+tv)$. We can also assume
that $|\gamma(1)|>R_0$ and all branching locus (other that $\infty$)
of $F_l(\xi,t):= \wh{F}_l(\xi, t v+X_\xi), G_j(\xi,t):=
\wh{G}_j(\xi, t v+X_\xi)$ are inside the disk $|t|<R_0$. Also, we
can assume that $\gamma$ does not hit the branching locus of
$F_l(\xi,t)$ and $G_j(\xi,t)$.

Write the Puiseux expansion of $h$ at $t = \infty$ as follows:

\begin{equation}\label{albert223}
h(\xi, t) = a_{\xi i_0} t^{\frac{i_0}{N_\xi}} (1+ o(1)), 
\end{equation}
where $N_\xi,\ i_0 \in {\mathbb N}, a_{\xi i_0} \not= 0$. 
We will fix $\xi \in U_0\sm E$ and $v$ as above for the rest of the
proof. Restricting (\ref{albert111}) to $X = t v+X_\xi$ and denoting
$\wh{F}_l(\xi, t v+X_\xi), \wh{G}_j(\xi, t v+X_\xi)$ by
$\wh{F}_l(\xi, t), \wh{G}_j(\xi, t)$ respectively for all $l, j$ as
before, one has the following equation  for $t\approx 0$:

\begin{equation}\label{albert222}
\prod_{l=1}^m \left(1+ \overline{F_l(z)} \cdot \wh{F}_l(\xi,
t)\right)^{\mu_l} = \-\phi(\bar z , \xi) \prod_{j=1}^v \left(1+
\overline{G_j(z)} \cdot \wh{G}_j(\xi, t)\right)^{\lambda_j} .
\end{equation}

If either a certain $\wh{F}_l(\xi, t)$ for some $l \not= 1$ or a
certain $\wh{G}_j(\xi, t)$ for some $j$, obtained by continuing
along $\gamma$ to $\gamma(1)$, has bounded branches by continuing in
the annulus  $R_0<|t|<\infty$, then there is an $\epsilon_0$ such
that for $|z|<\epsilon_0$, one has the following for $|t|>R_0$:
$$\frac{1}{2}< \left|1+ \overline{F_l(z)} \cdot \wh{F}_l(\xi, t)\right| \leq 2 ~~\text{or}~~\frac{1}{2} <
\left|1+ \overline{G_j(z)} \cdot \wh{G}_j(\xi, t)\right| \leq 2.$$
Otherwise, $$1+ \overline{F_l(z)} \cdot \wh{F}_l(\xi, t)
=  \sum_{\{i|s_i=s_l\}} \overline{f_{li}(z)} a_{\xi l s_i} t^{\frac{s_i}{N_\xi}} + o(t^{\frac{s_l}{N_\xi}}) 
~~\text{for}~t~\text{near}~\infty$$
 and for $s_i \in
\mathbb{Z}^{+}_0$, $s_l=\max_{i}s_i>0$, where $F_l(z)=(\cdots, f_{li}(z), \cdots)$,
$s_i/N_\xi$ is the order of blowing up of $f_{li}$. In addition, $\sum_{\{i|s_i=s_l\}}
\overline{f_{li}(z)} a_{\xi l s_i} \not \equiv 0$
%
by the arrangement that components in $F_l(z)$ are not linearly
dependent. We have a similar analysis for $1+ \overline{G_j(z)}
\cdot \wh{G}_j(\xi, t)$. Hence, for  unbounded branches, we have at
$t=\infty$:

\begin{equation}\label{new-002}
1+ \overline{F_l(z)} \cdot \wh{F}_l(\xi, t) = \-{a_{\xi l}(z)}
t^{\frac{s_l}{N_\xi}} + o(t^{\frac{s_l}{N_\xi}}),\ ~ 1+
 \overline{G_j(z)} \cdot \wh{G}_j(\xi, t) = \-{b_{\xi j}(z)} t^{\frac{t_j}{N_\xi}} + o(t^{\frac{t_j}{N_\xi}})
 \end{equation}
for $a_{\xi l}(z), b_{\xi j}(z) \not \equiv 0$ and $N_\xi, s_l, t_j
\in \mathbb{ N}$. Notice $\wh{F}_1$ is one of them with $s_1 >0$.
Now, we let $|z_0|<\epsilon_0$ be such that $a_{\xi l}(z_0), b_{\xi
j}(z_0) \not = 0$. Fix such a $z_0$. We perturb slightly $\gamma$
also to   make sure that  no terms in both sides of
(\ref{albert222}) hit zero. Moreover, the perturbed $\gamma$ has the
same terminal points and, relative to the terminal points, it has
the same homopotic class as the previous one in the space:
$${\CC}\sm \{ \hbox{branching points of } \wh{F}_l(\xi, t),
\wh{G_j}(\xi,t)\ \hbox{and}\ \ \hbox{zeros of }  1+
\overline{F_l(z_0)} \cdot \wh{F}_l(\xi, t),\ 1+
\overline{G_j(z_0)}\cdot \wh{G}_j(\xi, t) \ \hbox{in}\ t\}.$$

By holomorphic continuation of (\ref{albert222}) with $z=z_0$ along
the curve $\gamma(\tau)$ in $\tau$ from $ \tau=0$ to $\tau=1$, the
equation (\ref{albert222}) holds for certain fixed branches of
$\wh{F}_l(\xi, t)$ and $\wh{G}_j(\xi, t)$ with $z=z_0$  and $t$ near
$p_0=\gamma(1)$

  Hence, it follows
from equation (\ref{albert222}) with $z=z_0$, $t\approx p_0$, the
Puiseux expansions in (\ref{new-002}) centered at $t=\infty$, and
the blowing up rate at $t=\infty$ that

$$\prod_{S_l>0} \left( \-{a_{\xi l}(z_0)} t^{\frac{
s_l}{N_\xi}}\right)^{\mu_l} =c(z_0,\xi) \overline\phi(\bar z_0, \xi)
\prod_{t_j>0} \left( \-{b_{\xi j}(z_0)} t^{\frac{
t_j}{N_\xi}}\right)^{v_j}.$$
Here $s_1>0$ and  $s_l=0$, $t_j=0$ for those bounded branches. Also
$c(z_0,\xi)$ has the property that $$|c(z_0,\xi)|, \
\frac{1}{|c(z_0,\xi)|}<\infty.$$
Comparing the blowing-up rate for $t$ near $\infty$, one has:
$$\sum_{s_l>0} \mu_l s_l=\sum_{t_j>0} v_j t_j.$$
This contradicts the assumption in (\ref{james-10}).

 We thus proved  that all $\wh{F}_l(\xi,
X)$ and $\wh{G}_j(\xi, X)$ can not depend on $X$. Hence $F_l(\xi),
G_j(\xi)$ are algebraic. $\endpf$

\section{Single-valuedness of algebraic functions satisfying transcendental functional equations}

In this section, we derive the single-valuedness for algebraic
functions satisfying the transandental functional equation (\ref{log}).

\subsection{Starting point of the induction}

We first start with the following lemma:
\begin{lemma}\label{positive case}
Let $F_l: B \rightarrow \mathbb{C}^{N_l}$ be  algebraic and
holomorphic  for $1 \leq l \leq m$ with $B$ a ball centered at the
origin. Let  $h(z,\xi)$ be a rational and holomorphic  function in
$B$ such that $h(z,\-{z})>0$ for $z\in B$.
Suppose that
\begin{equation}\label{huang}
\prod_{l=1}^m \left(1+|F_l(z)|^2\right)^{\mu_l} = h(z,\-{z})
~~\text{on}~~B,
\end{equation}
for $\mu_l \in \mathbb{R}\sm \{0\}$. Then $F_l$ is a rational function
for each $1 \leq l \leq m$.
\end{lemma}

{\it Proof of Lemma \ref{positive case}:} There is a proper complex
analytic variety $E$ such that any branch of $F_l$ extends
holomorphically along any path $\gamma \subset {\mathbb C}^n\sm E.$
Moreover, for any $q\in {\mathbb C}^n\sm E$, the total number
$\nu(q)$ of all possible (different) values obtained by continuing
$F$ along closed curves in ${\mathbb C}^n\sm E$ starting from a
certain branch of $F$ near $q$ and then coming back and evaluating
at $q$, is independent of $q\in {\mathbb C}^n\sm E$.
 When
$\nu(q)=1$, $F$ is single-valued. Assume that $q\in B\sm E$ is
sufficiently close to $0$.
Write $U_l$ for an $(N_l+1)\times (N_l+1)$ unitary matrix such that
$$\wt{F_l}=(\wt{F}_{l,0},\cdots, \wt{F}_{l,N_l})=(1,F_l)\cdot U_l$$
has the property that $\wt{F_l}(q)=(c_l,0,\cdots, 0)$. Then
$|\wt{F}|^2=1+|F_l|^2$.
Write
$$G_l=(G_{l,1},\cdots,G_{l,N_l}):=(\frac{\wt{F}_{l,2}}{\wt{F}_{l,1}},\cdots,\frac{\wt{F}_{l,N_l}}{\wt{F}_{l,1}})$$
Write $h(z,\-{z})=h(q,\-{q})+2\Re{h_1(z)}+h_{mix}(z,\-{z})$ with
$h_1(z)=h(z,\bar q)-h(q,\-{q})$. Then the rational function $h_1$ takes
the zero value at $q$ and holomorphic in a neighborhood of $q$.
Moreover, the real-valued rational function $h_{mix}(z,\-{z})$ only
has mixed terms in its Taylor expansion at $(q,\-{q})$. Notice that
$h(q,\-{q})>0$ and $h_{mix}(q,\-{q})=0$. Taking
$\partial\-\partial\log$ to (\ref{huang}) and noticing that
$|\wt{F}|^2=1+|F_l|^2$, we get

\begin{equation}\label{huang-03}
\sum_{l=1}^m \mu_l\partial\-\partial\log
\left(1+|G_l(z)|^2\right)=\partial\-\partial\log{h}(z,\-{z}).
\end{equation}
Write

$$h(z,\-{z})=h(q,\-{q})(1+\frac{h_1(z)}{h(q,\-{q})})(1+\frac{\-{h_1(z)}}{h(q,\-{q})})\wt{h}(z,\-{z}).$$
Then $\wt{h}(z,\-{z})-1$ has only mixed terms in its Taylor
expansion at $(q,\-{q})$.
Since $\partial\-\partial\log
[(1+\frac{h_1(z)}{h(q,\-{q})})(1+\frac{\-{h_1(z)}}{h(q,\-{q})})]=0,$
comparing the Taylor expansion at $(q,\-q)$ for the left and right
hand sides of (\ref{huang-03}), after taking away
$\partial\-\partial$, we get
\begin{equation}\label{huang-04}
\prod_{l=1}^m \left(1+|G_l(z)|^2\right)^{\mu_l}=\wt{h}(z,\-{z}).
\end{equation}
Notice that $\wt{h}(q,\-{q})=1$. Now for any closed (simple) curve
$\gamma\in {\mathbb C}^n\sm E$ with $\gamma(0)=\gamma(1)=q,$ after
perturbing $\gamma$ slightly in the part $\gamma[\delta,1-\delta]$
for a small $\delta$ to avoid the zeros of $F_{l,0}$,
 $G_l(z)$ extends  holomorphically along $\gamma$. From the right
 hand side of (\ref{huang-04}),
we see that $\wt{h}(z,\-{z})$ is real analytic along $\gamma$, too.
 By the uniqueness
of real analytic functions, we see that (\ref{huang-03}) holds along
$\gamma$. Since $\wt{h}$ is single-valued and  $\wt{h}(q,\-{q})=1$,
we see that the value of $G_l$ continued along $\gamma$ still takes
$0$-value at $q$ for each $l$. (This is the only way to make sure
the right hand side of (\ref{huang-04}) attains its minimum value
$1$.) This shows that $G_l$ is single-valued for each $l$. Now,
notice that $ \wt{F_l}\cdot U_l^{-1}=(1,F_l)$. We have a certain
non-zero linear combination of ${\wt F_l}$ which is 1. Dividing on
both sides of such a linear combinantion by $F_{l,0}$, we see that
$F_{l,0}$ and thus $F_l$ are all single valued. Hence $F_l$ are rational. $\endpf$

\medskip

\begin{remark}
\label{jenny-002} There is an elegant argument by Mok in [M3] for
dealing with the single-valuedness of multi-valued maps, based on
the Cauchy-Schwarz inequality, which can also be used in our lemma
here.
Our argument for proving the single-valuedness of $F$ is based on
the fact that $1+|\cdot|^2$ has the only extremal value at $0$ and
the group $SU(n,1)$ acts transitively  on the quadratic form:
$|z|^2=|t|^2$. This argument applies to other type of potential
functions which may not have the Schwarz inequality property but
have similar minimum or maximum properties and symmetry. For
instance, it applies similarly to the following functional equation:
 \begin{equation}\label{func iden proj001}
 \prod_{l=1}^m \left(1-|F_l(z)|^2\right)^{\mu_l}=\chi(z,\-{z}), \
 \mu_l>0,\ \ |z|<<1.
\end{equation}
Here $F_l's$ are  holomorphic in a small neighborhood of $0$ with
$F_l(0)=0$, and extend along curves without hitting a certain proper
complex analytic variety, $\chi(z,\-{z})>0$ is  real analytic in
$|z|<1$ and $\chi(z,\xi)$ is meromorphic in $\{|z|<1\}\times
\{|\xi|<1\}$. Then one can also conclude the single-valuedness of
each $F_l$ in the unit ball.

\bigskip

\end{remark}

\begin{corollary}\label{proj}
Let $F_l: B \rightarrow \mathbb{C}^{N_l}$ be  a non-constant
algebraic and holomorphic map  for each $1 \leq l \leq m$ with a
ball $B$ centered at the origin. Suppose that there exists a
non-constant holomorphic irreducible polynomial function $h(z, \xi)$
over $B\times B$  with $h(z,\-{z})>0$ for $z\in {\CC}^n$
 such that

\begin{equation}\label{1111}
\prod_{l=1}^m \left(1+ F_l(z) \cdot \- {F_l(z)} \right)^{\mu_l} =
h(z, \-{z})^{r}\ \ \h{for} \ z\in B,
\end{equation}
and $\mu_l \in \mathbb{R}^+$,  ${r} \in \mathbb{R}$. There are a
positive integer $m_l$ and a positive constant $A_l$ such that $1+
F_l(z) \cdot \- {F_l(z)}=A_lh^{m_l}(z,\-{z})$ for each $l$.
\end{corollary}

{\it Proof of Corollary \ref{proj}:} After taking an
$r^{\rm{th}}$-root in (\ref{1111}), it follows from Lemma
\ref{positive case} that each $F_l$ is  a rational function, for $1
\leq l \leq m$. Since $h$ is a polynomial,
 $F_l$ is bounded on any compact subset of $\mathbb{C}^n$, $F_l$ thus is a
 polynomial.
Since the zeros of $(1+F_l(z)\bar{F}_l(\xi))$ can only take along
the zero set defined by $h(z,\xi)=0$, we finish the proof. $\endpf$

\bigskip

\subsection{Proof of Theorem \ref{singleness123}}

We now are in a position to prove the main result of this section:

\begin{theorem}\label{singleness123}
 Let $h(z, \xi)$ be
an irreducible polynomial   function over $\mathbb{C}^n \times
\mathbb{C}^n$ for $n \geq 1$ such that $h(z,\-{z})>0$.
Let $F_l(z): B\subset \mathbb{C}^n \rightarrow \mathbb{C}^{N_l},
G_j(z): B\subset\mathbb{C}^n \rightarrow \mathbb{C}^{N'_j}$ be
nonconstant holomorphic and algebraic over a small ball in ${\CC}^n$
centered at $0$ for $1 \leq l \leq m, 1 \leq j \leq v$.
Suppose the following transcendental functional identity holds:
\begin{equation}\label{trans identity0}
h(z, \-{z})^{{r}} = \prod_{l=1}^m \left( 1+ |F_l(z)|^2
\right)^{\mu_l} \cdot \prod_{j=1}^v \left( 1+ |G_j(z)|^2
 \right)^{-\lambda_j} ~~\text{over}~~B,\ \h{or~equivalently,}
\end{equation}
\begin{equation}\label{trans identity}
h(z, \xi)^{{r}} = \prod_{l=1}^m \left( 1+ F_l(z) \cdot \bar F_l(\xi)
\right)^{\mu_l} \cdot \prod_{j=1}^v \left( 1+ G_j(z) \cdot \bar
G_j(\xi)
 \right)^{-\lambda_j} ~~\text{over}~~B\times B,
\end{equation}
for  $\mu_l, \lambda_j \in \mathbb{R}^+$ satisfying
(\ref{james-10}), and for $r\in \RR\sm \{0\}$. Then $F_l, G_j$ are
holomorphic polynomials for all $l, j$. In fact, there exist
${m}_l, {n}_j \in \mathbb{N}$, such that
$$1+ F_l(z) \cdot \bar F_l(\xi) = A_l h(z, \xi)^{{m}_l} ~~\text{and}~ ~1+ G_j(z) \cdot \bar G_j(\xi) = B_j h(z, \xi)^{{n}_j},$$
for certain $A_l, B_j \in \mathbb{R}^+$. Moreover, $$r= \sum_{l=1}^m
{m}_l \mu_l - \sum_{j=1}^v \lambda_j {n}_j.$$
\end{theorem}

\medskip

Write $H_l(z, \xi) = 1+ F_l(z) \cdot \bar F_l(\xi)$ and
$Q_j(z,\xi)=1+ G_j(z) \cdot \bar G_j(\xi)$. Then $H_l(z, \xi),
Q_j(z,\xi)$ are holomorphic and algebraic  over $B\times B$.
Complexifying (\ref{trans identity0}), we get, after possibly
shrinking $B$ a little, the following:

\begin{equation}\label{isometry001}
h(z, \xi)^{r} = \prod_{l=1}^m H_l^{\mu_l}(z, \xi) \cdot
\prod_{j=1}^v Q_j^{-\lambda_j}(z, \xi)  \ \hbox{over } B\times B.
\end{equation}
Now, let $X_z$ be
the union of branching varieties of $F_l(z), G_j(z)$ and let $X$ be
the union of branching varieties of $H_l(z, \xi)$ and $G_j(z, \xi)$
for $l=1,\cdots,m, \ j=1,\cdots, v$. Then $F$ and $G$ can be
continued holomorphically along any path in $\mathbb{C}^n \setminus
X_z$.



\medskip

Our goal is to show that both $F$ and $G$ are forced to be
single-valued from the identity (\ref{isometry001}). Notice that the
argument in the proof of Corollary \ref{positive case} can not be
applied here, for now the potential functions involved
do not have any required extremal property. Indeed, without the
assumption in (\ref{james-10}), one can easily see that $F, G$ do
not have to be rational.


\begin{remark} We mention that there are many non-single-valued
 algebraic functions which only have zeros or points of indeterminancy  along the zero
set of an irreducible polynomial, say, $1+z \cdot \xi$. Indeed, let
$Y_1$ and $Y_2$  be the solutions of the polynomial $
Y^2+g(z)Y+(1+z\cdot \xi)$ with $g(0)=0$. Then $Y_1, Y_2$ are
holomorphic near $0$ and $Y_1\cdot Y_2=1+z\cdot \xi$. Hence, $Y_1,
Y_2$ can only have zeros in $1+z\cdot\xi=0$. Also, we can find
algebraic functions that have no zero at all, by simply considering
the solutions of polynomial equation $Y^2+(z\cdot\xi)Y+1=0$. Hence,
just by studying the distribution of zeros or points of
indeterminancy, we can not conclude the single-valuedness of
multi-valued maps.
\end{remark}

The proof of Theorem \ref{singleness123} is done by a long and
tedious induction argument based essentially on the monodromy
analysis. We will show that (\ref{james-10}) will be violated if
some $F_l$ or $G_j$ fails to be single-valued.

\medskip

{\it Proof of Theorem \ref{singleness123}}: We first notice that
under the assumption of Theorem \ref{singleness123}, all components
from $F_l$ and $G_j$ have leading coefficient $1$ in their minimal
polynomials.  We next fix more notations. When there are
$\alpha$-$H$'s and $\beta$-$Q$'s, we say that we are in the
situation with $(\alpha, \beta)$-factors. We use  ${r}_{*},
{r}_{**}$ to represent rational numbers, $n_{**}, m_{**}, k_{***}$
to denote integers, $n^+_{**}, m^+_{**}, k^+_{***}$ to denote
positive integers and $\iota_*, A$ to denote real numbers. All of
them may be different in different contexts.

\medskip

{\it Step 0: $H_l^{N_0}$ or $Q_j^{N_0}$ is a rational function.}
Assume $H_l(z, \xi)$ (resp. $Q_j(z, \xi)$) is a rational function
for some $l$ (resp. for some $j$). As all branches of $H_l$ (resp.
$G_j$) remain uniformly bounded over any
compact subset, then $H_l$ (resp. $G_j$) must be a polynomial. 
By Lemma \ref{leading coeff}, the zero  of $H_l(z, \xi)$ ( resp.
$Q_j(z, \xi)$) is contained in the variety defined by $\{(z, \xi):
h(z, \xi)=0\}$ (resp. $\{(z, \xi): h(z, \xi)=0\}$).
Hence, we conclude that $H_l(z, \xi) =A_l h(z, \xi)^{{m}^+_l}$ (resp. $Q_j(z, \xi)= B_j h(z, \xi)^{{n}^+_j}$). 

If for some $l$, all branching varieties of $H_l(z, \xi)$ are simple
cyclic branching varieties with respect to any branch, then there
exists an $N_0 \in \mathbb{N}$ such that $H_l^{N_0}$  is a
polynomial. Again, we get $H_l(z, \xi)^{N_0}
 =A'_l h(z, \xi)^{{m}^+_l}~$. Hence, $F_l$ is a polynomial by
 Corollary
 \ref{proj}.
 Now, if for each $l$, $H_l$ is single-valued,
 applying Corollary \ref{proj} again, we also conclude that $Q_j's$ and thus $G_j's$ are all
 single-valued polynomails. Hence the theorem is done in this setting. Similar
 arguments apply when for each $j$, $Q_j(z, \xi)$) has only simple cyclic branching varieties with respect to
any branch. We then apply
Lemma \ref{proj} to $H_l$ (resp. $Q_j$) to be single-valued. 

From now on, we always assume that $m,v\ge 1$, and also assume that
some $H_l$ and some $Q_j$ have non-cyclic branching varieties. We
will prove Theorem \ref{singleness123} by an induction argument on
$m$ and $v$. When $m=0$, we are done by Corollary \ref{positive
case}. We first consider the case of $(m, v)$-factor with $m=1, v
>0$ to illustrate our general argument.

\medskip

{\it Step 1: $(1, v)$-factor.}  We assume that Theorem
\ref{singleness123} is proved for the $(1, j)$-factor case with $0
\le j \le \beta-1$ and we consider the case with $(1,
\beta)$-factors. The following identity is the basic assumption:

\begin{equation}\label{calabi-001}
h^{{r}} = H_1^{\mu_1} Q_1^{-\lambda_1}\cdots
Q_\beta^{-\lambda_\beta} ~~~\text{near}~~~ 0.
\end{equation}
%


{\it Step 1.1:} 
As before, we can assume $H_1$ has an irreducible non-cyclic
branching hypervariety $E$
with respect to a certain branch of $H_1$. 
Applying Lemma \ref{branch} to (\ref{calabi-001}), we have

\begin{equation}\label{calabi-2}
m_{11}^+\mu_1=n_{11} \lambda_1+\cdots +n_{1\beta}\lambda_\beta,
\end{equation}
and
 \begin{equation}\label{jim16}
h^{m^+_{11} {r}} = (H_1^{-n_{11}}Q_1^{m^+_{11}})^{-\lambda_1} \cdots
(H_1^{-n_{1v}}Q_\beta^{m^+_{11}})^{-\lambda_\beta}.
\end{equation}

{\it Case 1.1.1:} If $(H_1^{-n_{11}}Q_1^{m^+_{11}})^{N_0}$ is
rational for some $N_0 \in \mathbb{N}$, then we have

\begin{equation}\label{solve-1}
H_1^{-n_{11}}Q_1^{m^+_{11}} =A h(z, \xi)^{{r}_1}
\end{equation}
for ${r}_1 \in \mathbb{Q}$, since the zeros and points of
indeterminancy of $(H_1^{-n_{11}}Q_1^{m^+_{11}})^{N_0}$ have to be
along $E_0= \{(z, \xi): h(z, \xi)=0\}$.

Now, if $n_{11}<0$, 
by Corollary  \ref{proj} with  $(2,
0)$-factors, we conclude that $Q_1=B_1 h(z, \xi)^{{r}_2}$ with ${r}_2 \in  \mathbb{N}$. 
Hence, after a cancelation, we reduce (\ref{calabi-001}) to the $(1,
\beta-1)$-factor case. Thus we are done by induction.

If $n_{11}>0$, then we get by (\ref{solve-1}):

\begin{equation}\notag
H_1= A Q_1^{\frac{m^+_{11}}{n_{11}}}h(z, \xi)^{{r}_2}
~~\text{or}~~Q_1=A H_1^{\frac{n_{11}}{m^+_{11}}} h(z, \xi)^{{r}_3}.
\end{equation}
Substituting into (\ref{calabi-001}), we get

\begin{equation}\notag
Q_1^{-\big(\lambda_1-\frac{m^+_{11}}{n_{11}}\mu_1\big)}Q_2^{-\lambda_2}
\cdots Q_\beta^{-\lambda_\beta}= A h(z, \xi)^{{\iota}_1}
~~\text{or}~~H_1^{\mu_1-\frac{n_{11}}{m^+_{11}}\lambda_1}
Q_2^{-\lambda_2} \cdots Q_\beta^{-\lambda_\beta}=A h(z,
\xi)^{{\iota}_2}.
\end{equation}
Furthermore, if $\frac{\lambda_1}{\mu_1} \geq
\frac{m^+_{11}}{n_{11}}$, then by Corollary \ref{proj}, we finish
the proof of Theorem \ref{singleness123} by inducntion. Otherwise,
we have $\frac{\lambda_1}{\mu_1}<\frac{m^+_{11}}{n_{11}}$, i.e.
$\mu_1
> \frac{n_{11}}{m^+_{11}} \lambda_1$. 
If $span_{\mathbb Q^+_0}
\bigg\{\mu_1-\frac{n_{11}}{m^+_{11}}\lambda_1\bigg\}
\cap span_{\mathbb Q^+_0} \left\{ \lambda_2, \cdots,
\lambda_\beta \right\}$=\{0\}, then we are also done
by the induction, as it is reduced to the $(1, \beta-1)$-factor
case. Otherwise, there exist nonnegative rational numbers $d^+_1,
c^+_2, \cdots, c^+_\beta$, not all zero, such that

\begin{equation}\notag
d_1^+\big(\mu_1-\frac{n_{11}}{m^+_{11}}\lambda_1\big)=\sum_{j=2}^\beta
c_j^+ \lambda_j > 0.
\end{equation}
Hence

\begin{equation}\notag
d_1^+ \mu_1 = d_1^+ \frac{n_{11}}{m^+_{11}} \lambda_1 +
\sum_{j=2}^\beta c_j^+ \lambda_1 \not=0.
\end{equation}
This contradicts (\ref{james-10}). This shows that $(H_1^{-n_{11}}Q_1^{m^+_{11}})^{N_0}$ is not a rational function for any $N_0 \in \mathbb{N}.$

\medskip

{\it Case 1.1.2: From $H_1^*Q_1^*$ to $H_1^*Q_1^*Q_2^*$.} Since $H_1^{-n_{11}}Q_1^{m^+_{11}}$ now has a non-cyclic branching variety, 
applying Lemma \ref{branch} to (\ref{jim16}), we have

\begin{equation}\notag
n^+_{11} \lambda_1=\sum_{j=2}^\beta n_{1j}\lambda_j, 
\end{equation}
and 

\begin{equation}\notag
h(z, \xi)^{m^+_{11}n^+_{11} {r}} = \bigg(H_1^{k_{22(-1)}}
Q_1^{k_{221}} Q_{2}^{k^+_{222}}\bigg)^{-\lambda_1} \cdots
\bigg(H_1^{k_{2\beta(-1)}} Q_1^{k_{2\beta1}}
Q_{\beta}^{k_{222}^+}\bigg)^{-\lambda_{\beta}}.
\end{equation}


{\it Step 1.2: From $H_1^*Q_1^*\cdots Q^*_p$ to $H_1^*Q_1^*\cdots Q^*_{p+1}$.} Now, 
suppose for some integer $p$ with $2 \leq p \leq \beta-1$, we have

\begin{equation}\label{comb6}
n^+_{(p-1)(p-1)} \lambda_{p-1} = \sum_{j=p}^\beta n_{(p-1)j} \lambda_j,
\end{equation}

and

\begin{equation}\label{jim1}
\begin{split}
h(z, \xi)^{m^+_{11}n^+_{11}\cdots n^+_{(p-1)(p-1)} {r}} &= \bigg(H_1^{k_{pp(-1)}} Q_1^{k_{pp1}}\cdots Q_{p-1}^{k_{pp(p-1)}} Q_{p}^{k_{ppp}^+}\bigg)^{-\lambda_{p}} \\
& \quad \cdots \bigg(H_1^{k_{p\beta(-1)}} Q_1^{k_{p\beta1}}\cdots
Q_{p-1}^{k_{pv(p-1)}} Q_{\beta}^{k_{ppp}^+}\bigg)^{-\lambda_{\beta}}
\end{split}
\end{equation}

\medskip

{\it Case 1.2.1:} Suppose

\begin{equation}\label{jim2}
H_1^{k_{pp(-1)}} Q_1^{k_{pp1}}\cdots Q_{p-1}^{k_{pp(p-1)}}
Q_{p}^{k_{ppp}^+} = A h(z, \xi)^{{r}_4}.
\end{equation}
We can assume $k_{pp(-1)} \leq 0$ after taking $(-1)$-th power in
(\ref{jim2}) if necessary. Since $k^+_{ppp}$ may be
 changed to negative, we will write $k_{ppp}$ later for $k^+_{ppp}$. For  simplicity of the notation, we assume
 that $k_{pp1}, \cdots, k_{ppk} >0$ for some integer $1\le k \leq p$ and others are nonpositive. If there is no such $k$,
 it is easily reduced to   the situation  with $(0, p+1)$-factors and we are done by Corollary \ref{proj}.
  Otherwise, by Lemma \ref{combo}, we can find $j'$ from $\{-1, 1, \cdots, k\}$ such that

\begin{equation}\label{comb5}
\frac{\lambda_{j}}{\lambda_{j'}} \geq \bigg|
\frac{k_{ppj}}{k_{ppj'}} \bigg|  ~~\text{for}~~j=-1, 1, \cdots, k.
\end{equation}
Here we write $\lambda_{-1}$ for $\mu_1$.

\medskip

{\it Case 1.2.1.1:} Assume $j' \not= -1$ and for simplicity of
notation, assume $j'=1$. By (\ref{jim2}), one solves

\begin{equation}\notag
Q_1= A h(z, \xi)^{{r}_5} H_1^{-\frac{k_{pp(-1)}}{k_{pp1}}}
Q_2^{-\frac{k_{pp2}}{k_{pp1}}} \cdots
Q_{p}^{-\frac{k_{ppp}}{k_{pp1}}}.
\end{equation}
Substituting  
into (\ref{calabi-001}), we have

\begin{equation}\label{jim4}
H_1^{\mu_1+ \frac{k_{pp(-1)}}{k_{pp1}} \lambda_1}
Q_2^{-\big(\lambda_2-\frac{k_{pp2}}{k_{pp1}}\lambda_1\big)} \cdots
Q_{p}^{-\big(\lambda_{p}-\frac{k_{ppp}}{k_{pp1}} \lambda_1\big)}
Q_{p+1}^{-\lambda_{p+1}} \cdots Q_\beta^{-\lambda_\beta} =A h(z,
\xi)^{{\iota}_3}.
\end{equation}
By (\ref{comb5}), the exponent of $H_1$ is non-negative and all
other exponents of $Q$'s are non-positive (with some strictly
negative). 
Now if
$$span_{\mathbb Q^+_0}\bigg\{\mu_1+ \frac{k_{pp(-1)}}{k_{pp1}} \lambda_1 \bigg\}
 \cap span_{\mathbb Q^+_0}\bigg\{\lambda_2-\frac{k_{pp2}}{k_{pp1}}\lambda_1, \cdots, \lambda_{p}-
 \frac{k_{ppp}}{k_{pp1}} \lambda_1, \lambda_{p+1}, \cdots, \lambda_\beta \bigg\}=\{0\},$$
  we can apply the induction hypothesis to make the conclusion.
  Otherwise, there exist $$d^+, c_1, c_2, \cdots, c_{\beta} \in \mathbb{Z}_0^+\ \hbox{with}\  d^+>0,$$
  such that

\begin{equation}
d^+ \mu_1 = \sum_{j=2}^{\beta} c_j \lambda_j - c_1 \lambda_1>0.
\end{equation}
Here we use  $\mathbb{Z}_0^+$ to denote the set of non-negative
integers. Also, we notice that $c_1
> 0$, for otherwise, by (\ref{james-10}), we have $d^+=0$, which is
a contradiction. Therefore, one gets

\begin{equation}
c_1 \lambda_1 = - d^+ \mu_1 + \sum_{j=2}^{\beta} c_j \lambda_j.
\end{equation}
Substituting it into (\ref{calabi-001}), we get

\begin{equation}
h(z, \xi)^{c_1{r}} = (Q_1^{d^+} H_1^{c_1})^{\mu_1} (Q_1^{c_2}
Q_2^{c_1})^{-\lambda_2} \cdots (Q_1^{c_{\beta}}
Q_{\beta}^{c_1})^{-\lambda_\beta}.
\end{equation}
By Lemma \ref{new-lemma} and the induction hypothesis with $(1,
\beta-1)$-factors, it follows that $Q_1^{d^+} H_1^{c_1} = A
h^{{r}_6}$. Hence $H_1=A_1 h(z, \xi)^{{r}_7}$ and we finish the
proof by induction.

\medskip

{\it Case 1.2.1.2:} Assume $j'=-1$.  We then solve from (\ref{jim2})

\begin{equation}
H_1 = A Q_1^{-\frac{k_{pp1}}{k_{pp(-1)}}} \cdots
Q_{p}^{-\frac{k_{ppp}}{k_{pp(-1)}}} h(z, \xi)^{{r}_{8}}.
\end{equation}
Substituting it into (\ref{calabi-001}), we get

\begin{equation}\label{jim5}
h(z, \xi)^{{\iota}_{4}} = A Q_1^{-\big(\lambda_1 +
\frac{k_{pp1}}{k_{pp(-1)}} \mu_1\big)} \cdots
Q_{p}^{-\big(\lambda_p+ \frac{k_{ppp}}{k_{pp(-1)}} \mu_1\big)}
Q_{p+1}^{-\lambda_{p+1}} \cdots Q_\beta^{-\lambda_\beta}.
\end{equation}
Notice that all exponents are  non-positive by (\ref{comb5}) with at least one negative. 
Applying 
Corollary \ref{proj}, it follows that $Q_j = B_j h^{ n_j^+}.$
We are also  done by induction.

\medskip

{\it Case 1.2.2:} Suppose  $H_1^{k_{pp(-1)}} Q_1^{k_{pp1}}\cdots
Q_{p-1}^{k_{pp(p-1)}} Q_{p}^{k_{ppp}^+}$ does not have all cyclic
branching varieties. By applying Lemma \ref{branch} once more,
we get (\ref{comb6}), (\ref{jim1}) with $p$ being  replaced by
$p+1$.

\medskip

{\it Step 1.3:}
By repeating the argument in Step 1.2, we will either prove the
theorem for $(1, \beta)$-factor or get (\ref{comb6}) for each $p
\leq \beta$ and at the last stage, we have (\ref{jim1}) for
$p=\beta$, i.e.

\begin{equation}\notag
h(z, \xi)^{m^+_{11}n^+_{11}\cdots n^+_{(\beta-1)(\beta-1)}{r}} =
(H_1^{k_{\beta \beta(-1)}} Q_1^{k_{\beta\beta 1}} \cdots
Q_\beta^{k^+_{\beta\beta\beta}})^{-\lambda_\beta}.
\end{equation}
Hence $\frac{\lambda_\beta}{r} \in \mathbb{Q}$. Back to
(\ref{comb6}) for  $2 \leq p \leq \beta$ and (\ref{calabi-2}), it
follows that $\frac{\mu_1}{r}, \frac{\lambda_1}{r}, \cdots,
\frac{\lambda_\beta}{r}$ are all rational numbers. This is a
contradiction to (\ref{james-10}). This proves the $(1,
\beta)$-factor case.

\medskip

Hence we complete the proof for the  $(1, v)$-factor case for any
$v$.

\medskip

{\it Step 2:} We assume that we already proved Theorem
\ref{singleness123} for $(l, j)$-factors with $1 \leq l< \alpha,
1\leq j\le  \beta$. We consider the case  with $(\alpha,
\beta)$-factors.  We have the following basic assumption:

\begin{equation}\label{isometry0001}
h(z, \xi)^{{r}} = H_1^{\mu_1} \cdots H_\alpha^{\mu_\alpha}
Q_1^{-\lambda_1} \cdots Q_\beta^{-\lambda_\beta}.
\end{equation}
We can assume that $H_1$ has a non-cyclic branching variety.
Then by applying Lemma \ref{branch} to (\ref{isometry0001}),  
we have

\begin{equation}\notag
n_{11}^+\mu_1= \sum_{l=2}^\alpha n_{1l}\mu_l + \sum_{j=1}^\beta m_{1j} \lambda_j, 
\end{equation}

and

$$h(z, \xi)^{n^+_{11}{r}}= \prod_{l=2}^\alpha \bigg( H_1^{k_{2(-l)(-1)}}H_l^{k^+_{2(-l)(-l)}} \bigg)^{\mu_l} \cdot \prod_{j=1}^\beta \bigg( H_1^{k_{2j(-1)}} Q_j^{k^+_{2jj}} \bigg)^{-\lambda_j} .$$


{\it Step 2.1:} Suppose for $p$ with $2 \leq p \leq \alpha$, we have

\begin{equation}\label{jim9}
n^+_{(p-1)(p-1)} \mu_{p-1} = \sum_{l=p}^\alpha n_{(p-1)l}\mu_l + \sum_{j=1}^\beta m_{(p-1)j} \lambda_j, 
\end{equation}
and

\begin{equation}\label{jim10}
\begin{split}
h(z, \xi)^{n^+_{11}\cdots n^+_{(p-1)(p-1)}{r}} = &
\prod_{l=p}^{\alpha} \bigg(H_1^{k_{p(-l)(-1)}} \cdots
H_{p-1}^{k_{p(-l)(-p+1)}}
 H_{l}^{k^+_{p(-l)(-l)}}\bigg)^{\mu_l} \\
&\cdot \prod_{j=1}^{\beta} \bigg(H_1^{k_{pj (-1)}} \cdots
H_{p-1}^{k_{pj(-p+1)}} Q_{j}^{k^+_{pjj}}\bigg)^{-\lambda_j}.
\end{split}
\end{equation}


{\it Case 2.1.1:} Suppose that all branching varieties of
$H_1^{k_{p(-p)(-1)}} \cdots H_{p-1}^{k_{p(-p)(-p+1)}}
H_{p}^{k^+_{p(-p)(-p)}}$ are simple cyclic branching varieties with
respect to any of its branch, then

\begin{equation}\label{jim6}
H_1^{k_{p(-p)(-1)}} \cdots H_{p-1}^{k_{p(-p)(-p+1)}}
H_{p}^{k^+_{p(-p)(-p)}} =A h(z, \xi)^{{r}_{10}}.
\end{equation}
If $k_{p(-p)(-1)}, \cdots, k^+_{p(-p)(-p)} \geq 0$, then we are done
by applying
 Corollary \ref{proj}.
We thus assume that at least one exponent is negative. For
simplicity of notation, assume $k_{p(-p)(-1)}, \cdots, k_{p(-p)(-k)}\\
< 0$ and $k_{p(-p)(-k-1)}, \cdots, k_{p(-p)(-p)} \geq 0$ for $1 \leq
k \leq p-1$. By Lemma \ref{combo}, we can find $1 \leq l' \leq k$
such that

\begin{equation}
\frac{\mu_l}{\mu_{l'}} \geq \bigg|
\frac{k_{p(-p)(-l)}}{k_{p(-p)(-l')}} \bigg| ~~\text{for}~~ 1 \leq l
\leq k.
\end{equation}
%
We further assume $l'=1$ to simplify the notation. Then we solve
from (\ref{jim6})

\begin{equation}\notag
H_1=A H_2^{-\frac{k_{p(-p)(-2)}}{k_{p(-p)(-1)}}} \cdots
H_{p}^{-\frac{k_{p(-p)(-p)}}{k_{p(-p)(-1)}}} h(z, \xi)^{{r}_{11}}.
\end{equation}
Substituting it  into (\ref{isometry0001}), we have

\begin{equation}
H_2^{\mu_2-\frac{k_{p(-p)(-2)}}{k_{p(-p)(-1)}} \mu_1} \cdots
H_{p}^{\mu_{p}-\frac{k_{p(-p)(-p)}}{k_{p(-p)(-1)}}\mu_1}
H_{p+1}^{\mu_{p+1}} \cdots H_{\alpha}^{\mu_\alpha} Q_1^{-\lambda_1}
\cdots Q_\beta^{-\lambda_\beta} =A h(z, \xi)^{{\iota}_5}.
\end{equation}
Notice that all the exponents for $H_2, \cdots, H_\alpha$ are
nonnegative. 
Now either
we are in the case of

\begin{equation}\label{comb7}
\begin{split}
span & _{\mathbb Q^+_0}  \bigg\{ \mu_2-\frac{k_{p(-p)(-2)}}{k_{p(-p)(-1)}}\mu_1, \cdots, \mu_{p}-\frac{k_{p(-p)(-p)}}{k_{p(-p)(-1)}}\mu_1, \mu_{p+1}, \cdots, \mu_\alpha \bigg\} \\
&\cap span_{\mathbb Q^+_0} \big\{ \lambda_1, \cdots, \lambda_\beta
\big\} = \{ 0\}
\end{split}
\end{equation}
and we then finish the proof of the theorem by the induction
hypothesis, or (\ref{comb7}) does not hold. Therefore, there exist
$d_1, d_2, \cdots, d_\alpha, c_1, \cdots, c_\beta \in \mathbb
Q^+_0$, such that

\begin{equation}\notag
\sum_{l=2}^{\alpha} d_l \mu_l - d_1\mu_1 = \sum_{j=1}^\beta c_j
\lambda_j >0.
\end{equation}
It follows that $d_1 > 0$ by (\ref{james-10}). Hence, as before, we
get

\begin{equation}\label{comb8}
d_1 \mu_1 = \sum_{l=2}^\alpha d_l \mu_l - \sum_{j=1}^\beta c_j
\lambda_j.
\end{equation}
Taking $d_1$-th power in (\ref{isometry0001}) and using
(\ref{comb8}), it follows that

\begin{equation}\notag
h(z, \xi)^{d_1 {r}} = \prod_{l=2}^\alpha
(H_1^{d_l}H_l^{d_1})^{\mu_l} \cdot \prod_{j=1}^\beta
(H_1^{c_j}Q_j^{d_1})^{-\lambda_j}.
\end{equation}
By the same argument at the end of Case 1.2.1.1, we are done
applying 
the induction hypothesis.

\medskip

{\it Case 2.1.2:} Suppose (\ref{jim6}) does not hold. 
By 
Lemma \ref{branch}, we get (\ref{jim9}), (\ref{jim10}) for $p$ being
replaced by $p+1$. By repeating the same argument as in Step 2.1, we
will either prove the theorem, or get (\ref{jim9}) for each $2 \leq
p \leq \alpha+1$ and (\ref{jim10}) for $p=\alpha+1$. Hence,

\begin{equation}\label{jim100}
n^+_{\alpha\alpha} \mu_{\alpha} =  \sum_{j=1}^\beta m_{\alpha j} \lambda_j, 
\end{equation}
and

\begin{equation}\notag
h(z, \xi)^{n^+_{11}\cdots n^+_{\alpha \alpha} {r}} =
\prod_{j=1}^{\beta} \bigg(H_1^{k_{(\alpha+1) j (-1)}} \cdots
 H_{\alpha}^{k_{(\alpha+1) j(-\alpha)}} Q_{j}^{k^+_{(\alpha+1)jj}}\bigg)^{-\lambda_j}.
\end{equation}
By applying Lemma \ref{branch}, we either conclude
$$H_1^{k_{(\alpha+1) 1 (-1)}} \cdots H_{\alpha}^{k_{(\alpha+1) 1(-\alpha)}} Q_{1}^{k^+_{(\alpha+1)11}} =
A h(z, \xi)^{{r}_{12}}$$ or
$$n^+_{(\alpha+1)(\alpha+1)}\lambda_{1} = \sum_{l=2}^\beta m_{(\alpha+1)l} \lambda_l,$$ and

$$h(z, \xi)^{n^+_{11}\cdots n^+_{(\alpha+1)(\alpha+1)} {r}} = \prod_{j=2}^{\beta} \bigg(H_1^{k_{(\alpha+2) j (-1)}}
 \cdots H_{\alpha}^{k_{(\alpha+2) j(-\alpha)}} Q_{1}^{k_{(\alpha+2) j j}}  Q_{j}^{k^+_{(\alpha+2)jj}}\bigg)^{-\lambda_j}.$$


{\it Step 2.2:} Suppose we know for some integer $2 \leq p \leq
\beta-1$,


\begin{equation}\label{jim11}
n^+_{(\alpha+p-1)(\alpha+p-1)}\lambda_{p-1} = \sum_{l=p}^\beta m_{(\alpha+p-1)l} \lambda_l 
\end{equation}
and

\begin{equation}\label{jim14}
h(z, \xi)^{n^+_{11}\cdots n^+_{(\alpha+p-1)(\alpha+p-1)}{r}} =
\prod_{j=p}^{\beta} \bigg(H_1^{k_{(\alpha+p) j (-1)}}
 \cdots H_{\alpha}^{k_{(\alpha+p) j(-\alpha)}} Q_{1}^{k_{(\alpha+p) j 1}} \cdots Q_{j}^{k^+_{(\alpha+p)jj}}\bigg)^{-\lambda_j}.
\end{equation}
By applying Lemma \ref{branch}, we either conclude

\begin{equation}\label{jim12}
 H_1^{k_{(\alpha+p) p(-1)}} \cdots H_{\alpha}^{k_{(\alpha+p) p(-\alpha)}} Q_{1}^{k_{(\alpha+p)p 1}}
 \cdots Q_{p}^{k^+_{(\alpha+p)pp}}=A h(z, \xi)^{{r}_{13}}
\end{equation}
or (\ref{jim11}), (\ref{jim14}) with $p$ being replaced  by $p+1$.

\medskip

{\it Case 2.2.1:} Assume (\ref{jim12}) for $1 \leq p \leq \beta-1$.
If $k_{(\alpha+p)p (-1)}= \cdots = k_{(\alpha+p)p (-\alpha)}= 0$,
 then we are done by repeating exactly the same argument as in Case 2.1.1. Otherwise,
  by taking $(-1)$-th power of (\ref{jim12}) if necessary, we can assume $k_{(\alpha+p)
  p(-1)}<0$. If this is the case, $k^+_{(\alpha+p)pp}$ will change
  the sign and hence we will write in the following
  $k_{(\alpha+p)pp}$ for $k^+_{(\alpha+p)pp}$.
Suppose, without loss of generality, that there exists $1\leq k \leq
\alpha, 0 \leq k' \leq p$, such that $k_{(\alpha+p)p (-l)}<0$ for $1
\leq l \leq k$, $k_{(\alpha+p)p(-l)} \geq0$ for $l > k$, and
$k_{(\alpha+p)p j}>0$ for $1\leq j \leq k'$, $k_{(\alpha+p)pj} \leq
0$ for $p \geq j > k'$. Here $k' =0$ means that $k_{(\alpha+p)pj}
\leq 0$ for all $j=1, \cdots, p$. By Lemma \ref{combo}, there exists
$-1\geq j' \geq -k$ or $1 \leq j' \leq k'$ such that

\begin{equation}
\frac{\lambda_j}{\lambda_{j'}} \geq \bigg|
\frac{k_{(\alpha+p)pj}}{k_{(\alpha+p)pj'}} \bigg|~~\text{for}~~j=-1,
\cdots, -k ~\text{and}~ 1, \cdots, k',
\end{equation}
Here we use $\lambda_{j}$ to denote $\mu_{-j}$ when $j<0$.

\medskip

{\it Case 2.2.1.1:  $j'<0$}. Without loss of generality and for
simplicity of notation, assume $j'=-1$. Then we solve from
(\ref{jim12})

\begin{equation}\notag
H_1 = A H_2^{-\frac{k_{(\alpha+p)p(-2)}}{k_{(\alpha+p)p(-1)}}}
\cdots
H_\alpha^{-\frac{k_{(\alpha+p)p(-\alpha)}}{k_{(\alpha+p)p(-1)}}}
Q_1^{-\frac{k_{(\alpha+p)p1}}{k_{(\alpha+p)p(-1)}}} \cdots
Q_p^{-\frac{k_{(\alpha+p)pp}}{k_{(\alpha+p)p(-1)}}} h(z,
\xi)^{{r}_{14}}.
\end{equation}
Substituting it  into (\ref{isometry0001}), we have

\begin{equation}\notag
\prod^\alpha_{l=2} H_l^{\mu_l
-\frac{k_{(\alpha+p)p(-l)}}{k_{(\alpha+p)p(-1)}}\mu_1} \cdot
\prod_{j=1}^p
Q_j^{-\big(\lambda_j+\frac{k_{(\alpha+p)pj}}{k_{(\alpha+p)p(-1)}}\mu_1\big)}
\cdot \prod_{j=p+1}^\beta Q_{j}^{-\lambda_{j}} = A h(z,
\xi)^{{\iota}_{6}}.
\end{equation}
Notice that  the exponents for $H_2, \cdots, H_\alpha$ are
nonnegative and those of $Q_1, \cdots, Q_\beta$ are nonpositive and
at least one is not zero. 
Now either we are in the case of

\begin{equation}\label{comb9}
\begin{split}
span & _{\mathbb Q^+_0}  \bigg\{ \mu_2-\frac{k_{(\alpha+p)p(-2)}}{k_{(\alpha+p)p(-1)}}\mu_1, \cdots, \mu_\alpha-\frac{k_{(\alpha+p)p(-\alpha)}}{k_{(\alpha+p)p(-1)}}\mu_1 \bigg\} \\
&\cap span_{\mathbb Q^+_0}  \bigg\{
\lambda_1+\frac{k_{(\alpha+p)p1}}{k_{(\alpha+p)p(-1)}}\mu_1, \cdots,
\lambda_p+\frac{k_{(\alpha+p)pp}}{k_{(\alpha+p)p(-1)}}\mu_1,
\lambda_{p+1}, \cdots, \lambda_\beta \bigg\} = \{ 0\}
\end{split}
\end{equation}
and we finish the proof of the theorem by the induction hypothesis,
or (\ref{comb9}) does not hold. Therefore, there exist $d_1, d_2,
\cdots, d_\alpha, c_1, \cdots, c_\beta \in {\mathbb Z}^+_0$, such
that

\begin{equation}\notag
\sum_{l=2}^{\alpha} d_l \mu_l - d_1\mu_1 = \sum_{j=1}^\beta c_j
\lambda_j >0.
\end{equation}
It follows that $d_1 > 0$ by (\ref{james-10}) and we finish our
proof as in Case 2.1.1.

\medskip

{\it Case 2.2.1.2:} {\it  $j'>0$}. Without loss of generality,
assume $j'=1$. Repeat the same argument as in {\it Case 2.2.1.1}, we
also finish the proof.

\medskip

{\it Case 2.2.2:} Suppose we have (\ref{jim11}) (\ref{jim14}) with
$p$ being replaced  by $p+1$. By repeating the same argument as in
Case 2.2.1, we eventually have (\ref{jim11}) for all $2 \leq p \leq
\beta$ and (\ref{jim14}) for $p=\beta$. Hence,

\begin{equation}\notag
h(z, \xi)^{n^+_{11} \cdots
n^+_{(\alpha+\beta-1)(\alpha+\beta-1)}{r}} =
\big(H_1^{k_{(\alpha+\beta)\beta(-1)}} \cdots
H_1^{k_{(\alpha+\beta)\beta(-\alpha)}}
Q_1^{k_{(\alpha+\beta)\beta1}} \cdots
Q_\beta^{k^+_{(\alpha+\beta)\beta\beta}} \big) ^{-\lambda_\beta}.
\end{equation}
Hence we conclude $\frac{\lambda_\beta}{r} \in \mathbb{Q}$.

\medskip

{\it  Step 2.3:} Since (\ref{jim11}) holds for $2 \leq p \leq
\beta$, (\ref{jim100}) and (\ref{jim9}) hold for $2 \leq p \leq
\alpha$, it follows that $\frac{\mu_1}{r}, \cdots,
\frac{\mu_\alpha}{r},  \frac{\lambda_1}{r}, \cdots,
\frac{\lambda_\beta}{r} \in \mathbb{Q}$. We arrive at a
contradiction
 with (\ref{james-10}). This proves the $(\alpha, \beta)$-factor case.

By  induction, we finally complete the proof of Theorem
\ref{singleness123}. $\endpf$


\subsection{Proof of Theorem \ref{kahler}}

As an application of Theorem \ref{singleness123}, we give a proof of
Theorem \ref{kahler}.

\medskip {\it Proof of Theorem \ref{kahler}}:  Suppose the notation
and assumption in Theorem \ref{kahler}. $F_l, G_j$ can be regarded
as holomorphic maps from $U\subset {\CC}^n$ into ${\mathbb C}^N$
with a large $N$, after shrinking $U$ and choosing suitable
coordinates. We can assume that $F_l(0), G_j (0)=0$.  Then
(\ref{kahler isometry}) is equivalent to the following:

\begin{equation} \label{james-0007} \sqrt{-1}  \partial\bar\partial
\log h(z,\-{z}) =\sqrt{-1} \sum_{l=1}^m \mu_l
\partial \bar\partial \log(1+|F_l(z)|^2)-\sqrt{-1} \sum_{j=1}^v \lambda_j
\partial \bar\partial \log(1+|G_j(z)|^2),
 \ z\in U.
\end{equation}

Taking away $\p\-{\p}$ and comparing the mixed terms in the Taylor
expansion of (\ref{james-0007}) around $0$, we get (\ref{trans
identity0}). By Theorems \ref{algebraicity-1} and
\ref{singleness123}, we conclude that  we have $(1+|F_l|^2)=A_l
h^{m_l}$ and $(1+|G_j|^2)=B_j h^{n_j}$ for each $l,j$ with
$m_l,n_j\in {\mathbb N},\  A_l, B_j
>0$. Hence, $F_l: (U\subset M, m_j\omega)\ra ({\mathbb P}^N,
\omega_{FS})$ and $G_j: (U\subset M, n_j\omega)\ra ({\mathbb P}^N,
\omega_{FS})$ are local isometries for all $l, j$. By the classical
result of Calabi in [Ca], $F_l, G_j$ extend holomorphically along
any path inside $M$. Since $M$ is simply connected, we see $F_l$ or
$G_j$ extend to holomorphic immersions from $M$ into ${\mathbb
P}^{N_l}$ or ${\mathbb P}^{N_j'}$, respectively,  with the same kind
of isometric property.
The rest of the proof is easy. $\endpf$

\medskip

We mention that Theorem \ref{kahler} applies when $(M,\omega)$ is
${\mathbb P}^n$ with the standard Fubini-Study metric.  Indeed, as
even in the more general compact Hermitian symmetric spaces of
compact type case, the extended maps are one-to-one.
(See the proof of Theorem \ref{maintheorem-local} in the following
section.)


\section{Proof of Theorem \ref{maintheorem-local}}
Let $M\subset {\mathbb P}^n$ be a projective algebraic manifold of
complex dimension $k$. Let $\{U_j,\psi_j\}$ be a system of
holomorphic charts of $M$ with $\psi_j: U_j\ra V_j\subset {\mathbb
C}^k$ a biholomorphic map for each $j$. (We always assume that
$U_j's$ are connected.) We call $\{U_j,\psi_j\}$ is a system of Nash
algebraic holomorphic charts if for each $j$ and writing
$\psi_j^{-1}=[\phi_{j0},\cdots,\phi_{jn}]$ with
$\phi_{jl_0}\not\equiv 0$, then $\frac{\phi_{jl}}{\phi_{jl_0}}$ is
Nash algebraic over $V_j$ for each $l$.
A holomorphic chart of $M$ is  said to be a  Nash algebraic
holomorphic chart if the resulting transition functions with respect
to $\{U_j, \psi_j\}$ are Nash algebraic.  Now, a holomorphic
function $h$ defined over a connected open subset $U$ of $M$ is
called Nash-algebraic if for any $U_j$ with $U_j\cap U\not
=\emptyset,$ $h\circ\psi_j^{-1}$ is Nash algebraic over
$\psi_j(U_j\cap U)$. Suppose both $M\subset {\mathbb P}^n$ and
$M'\subset {\mathbb P}^{n'}$ are projective algebraic manifolds with
holomorphic Nash-algebraic systems $\{U_j,\psi_j\}$ and
$\{U'_j,\psi'_j\}$, respectively . A holomorphic map $G: U\subset
M\ra M'$ is called holomorphic Nash algebraic, if $\psi'_k\circ
G\circ (\psi_j)^{-1}$ is Nash algebraic whenever it is well defined.
Apparently, this definition is independent of the choice of the
systems of Nash algebraic charts. Also, the Nash-algebraicity is a
global property in the sense that if $G$ is Nash algebraic in a
small open subset of $U$, then it is Nash algebraic over $U$.

By a
variation of Chow's theorem ([p.170, GH]), when $G$ is a global map
from $M$ into $M'$, then G is a restriction of a rational map from
${\mathbb P}^{n}$ into ${\mathbb P}^{n'}$, and thus is Nash
algebraic in our definition. Also, all the above definitions are
independent of the embedding from the algebraic manifolds into
projective spaces.

 Let $(M,\omega_M)$ be an irreducible Hermitian symmetric space
of compact type of complex dimension $n$ with a fixed canonical
K\"ahler-Einstein metric $\omega_M$.
We also identify  $\omega_M$ with its associated
K\"ahler form if there is no risk of causing confusion. 
 By the Nakagawa-Takagi isometric embedding
theorem [M1], we can assume that $M$ has been holomorphically
isometrically embedded into a certain projective space $({\mathbb
P}^{n_0}, \omega_{n_0})$ equipped with the standard Fubini-Study
metric $\omega_{n_0}$. Let $U$ be a connected open subset of $M$.
Consider  local holomorphic  mappings $F=(F_1,\cdots,F_m):
U\rightarrow {\mathbb P}^{N_1}\times\cdots\times {\mathbb P}^{N_m}$
and $G=(G_1,\cdots,G_v): U\rightarrow {\mathbb
P}^{N'_1}\times\cdots\times {\mathbb P}^{N'_v}$ such that

\begin{equation}\label{james-01}
\omega_M =\sum_{l=1}^{m}\mu_l F^{*}_l
\omega_{N_l}-\sum_{j=1}^{v}\lambda_j G^{*}_j \omega_{N'_j},
\end{equation}
 where for each $j,l$, the numbers
$\lambda_j, \mu_l$ are  positive real numbers satisfying the number
theoretic property in (\ref{james-10}), and $\omega_{N_l}$,
$\omega_{N'_j}$ are the standard Fubini-Study metric over ${\mathbb
P}^{N_l}$ and ${\mathbb P}^{N_j'},$ respectively.

Since  $(M,\omega_M)\subset ({\mathbb P}^{n_0},\omega_{n_0})$ is a
projective algebraic manifold, it is defined by a set of homogeneous
polynomials. We can thus  find a finite set of holomorphic
coordinate charts  $\{U_i, \Psi_i\}$, that covers $M$. We  assume
that each $U_i$ is contained in a standard complex coordinate chart
of ${\mathbb P}^{n_0}$ and each $\Psi_i$ is obtained by  mapping
into ${\CC}^{n_0}$ in a standard way and then  projecting down to  a
certain open subset $V_i$ in the complex Euclidean space ${\mathbb
C}^n$ with coordinates $z=(z_1,\cdots,z_n)$. Write the inverse of
$\Psi_i$ to be $\Phi_i=[\phi_{i0},\cdots,\phi_{i{n_0}}]$. Hence,
there is a certain component  $\phi_{ij}\not = 0$ over $V_i$.
We can assume, without loss of generality, that each component of
$\Phi_i$ is holomorphic Nash algebraic.
The system of holomorphic coordinate charts thus obtained is a
Nash-algebraic system. Notice that the metric $\omega_M$ is
represented by
 $\sum_{\alpha,\beta=1}^{n}h_{\alpha\-{\beta}}dz_\alpha \otimes d\-{z}_\beta$ on $V_i$ with
 $h_{\alpha\-{\beta}}={\p^2 \over \p z_\alpha \p \-{z}_\beta}\log
 (\sum_{k=0}^{n_0}|\phi_{ik}|^2).$

Write
$\omega^{*}=\sqrt{-1}\partial\-{\partial}\log(\sum_{k=0}^{n_0}|\phi_{ik}|^2).$
Suppose, without loss of generality, that
 $U\subset U_1$.  Write $U^{*}=\Psi_1(U) \subset V_1$.
 Shrinking $U$ if necessary, we can  further  assume that
for each $l$, $F_l(U)$ is contained in one of the $N_l+1$ standard
holomorphic coordinate charts $V_{N_l, l_0}$ of ${\mathbb P}^{N_l}$
with the standard (linear-fractional) holomorphic coordinate map
$\sigma_{N_l,l_0}: V_{N_l,l_0}\ra {\mathbb C}^{N_l}$, and we can
assume that $G_j(U)$ is contained in one of the $N'_j+1$ standard
holomorphic coordinate charts $V_{N'_j, j'_0}$ of ${\mathbb
P}^{N'_j}$ with the standard  coordinate map $\sigma_{N'_j,j'_0}:
V_{N'_j,j'_0}\ra {\mathbb C}^{N'_j}$. Here, we recall that
$V_{N_j,k}=\{[z_0,\cdots,z_{N_j}]:\ z_k\not = 0\}$ and
$\sigma_{N_j,k}([z_0,\cdots,z_{N_j}])=(z_0/z_k,\cdots,z_{k-1}/z_k,z_{k+1}/z_k,\cdots,z_{N_j}/z_k).$

\medskip

{\it Proof of Theorem \ref{maintheorem-local}}: Assume the notation
and assumptions in Theorem \ref{maintheorem-local} and above.  As
mentioned before, by the isometric embedding theorem of
Nakagawa-Takagi ([Theorem 1, p.135, M1]), we need only consider maps
into the product of complex projective spaces.

We now let $M$ be an irreducible Hermitian symmetric space of
compact type of dimension $n$. As above, we assume that $M$ is
holomorphically isometrically embedded into $\mathbb{P}^{n_0}$ for a
certain $n_0 \geq 2$. We let $(U_1, \Psi_1)$ be a Nash algebraic
holomorphic coordinate chart of $M$ with holomorphic Nash algebraic
maps
$$\Phi_1 = \Psi^{-1}_1: V_1 = \Psi_1(U_1) \subset \mathbb{C}^n \rightarrow \mathbb{P}^{n_0}$$
such that $U \subset U_1$. We also assume that $U_1$ is  contained
in one of standard coordinate charts of ${\mathbb P}^{n_0}$, namely,
one of the coordinate in a certain fixed order never vanishes for
points in $U_1$.

Choose a point $p^*\in U$ and after composing with isometries, we
can assume $F_l(p^*)$ and $G_j(p^*)$ are all mapped by the standard
coordinate maps to the origin in the complex Euclidean spaces.
Moreover, after shrinking $U$, and thus $U^*:=\Psi_1(U)$, we can
assume $F_l(U^*) \subset V_{N_l,l_0}, G_j(U^*) \subset
V_{N'_j,j'_0}.$ Also, write
$$Y=(Y_1,\cdots,Y_{m+v}):=$$
$$(\sigma_{N_1, 1_0}\circ F_1 \circ \Psi_1^{-1}, \cdots, \sigma_{N_m,
m_0}\circ F_p \circ \Psi_1^{-1}, \sigma_{N'_1,1'_0}\circ G_1 \circ
\Psi_1^{-1},\cdots, \sigma_{N'_v, v'_0}\circ G_v \circ \Psi_1^{-1})
   \big|_{\Psi_1(U)},$$
 where  $\sigma_{N_j,j_0}$ and $  \sigma_{N'_j,j'_0}$ are standard coordinate maps.
Then  $Y(\Psi_1(p^*))=0$.

After applying a holomorphic isometry, we also assume that
$p^*=[1,0,\cdots,0]$. Shrinking $U_1$ and $U$, we can assume that
the first coordinate of points in $U_1$ is never zero. Write
$\Psi^{-1}_1=[1,\phi_1(z),\cdots, \phi_{n_0}(z)]$. Here, $\phi_k$ is
algebraic and  holomorphic over $V_1:=\Psi_1(U_1)$. Without loss of
generality, we  identify $p^*$  with the origin in ${\CC}^{n_0}$ and
$\phi_k(0)=0$ through the coordinate map $\Psi_1$.
From the
hypothesis (\ref{james-01}), it follows that
\begin{equation} \label{james-007}
\sqrt{-1}  \partial\bar\partial
\log(1+\sum_{k=1}^{n_0}|\phi_{k}(z)|^2)+\sqrt{-1} \sum_{j=1}^v
\lambda_j
\partial \bar\partial \log(1+|G_j(z)|^2) =\sqrt{-1} \sum_{l=1}^m \mu_l
\partial \bar\partial \log(1+|F_l(z)|^2).
\end{equation}
Here,
for simplicity of notation,  
we still write $F$ for $(\sigma_{N_1,1_0}\circ
F_1,\cdots,\sigma_{N_m,m_0}\circ F_m)\circ \Phi_1$ and $G$ for
$(\sigma_{N'_1,1'_0}\circ G_1,\cdots,\sigma_{N'_v,v'_0}\circ
G_v)\circ \Phi_1$. As in Clozel-Ullmo [CU] and Mok [M2],
(\ref{james-007}) yields the following:

\begin{equation}\label{func iden proj2-d}
\prod_{l=1}^m  \left(1+|Y_l|^2 \right)^{\mu_l}\cdot \prod_{j=1}^v
\left( 1+|Y_{m+j}|^2 \right)^{-\lambda_j}=1+\sum_{k}|\phi_k(z)|^2, \
\ z\in \Psi_1(U).
\end{equation}
 Now, by Theorem \ref{algebraicity-1}, we know that $Y$ is algebraic over $V_1$.

    We next prove that $Y$ admits a holomorphic extension along any path in $V_1$.
    Suppose not.
     We will have an irreducible branching variety $E$ for (some component of) $Y$ in $V_1$.
     Let $p_0^* \in E$ be a generic smooth point.
     Then after an algebraic holomorphic change of coordinates, we can assume that $p_0^*=0$ and $E$ near 0
      is defined by $z_n=0$. Now
      for $p^*$ (near 0) $\not\in E$, we have a certain branch $\wh{Y}$  of $Y$ near $p^*$,
       which extends to a multi-valued holomorphic map in a small neighborhood of $0$ and admits the following Puiseux  expansion:

$$\wh{Y} = \sum_{\alpha=0}^\infty \tilde{a}_\alpha(z') z_n^{\frac{\alpha}{N}}.$$
Here we can find the smallest positive integer $\alpha_0$ with $N
\nmid \alpha_0$ such that the vector-valued
$\tilde{a}_{\alpha_0}(z') \not\equiv 0$ for $z'$ near 0. By the
uniqueness of real analytic functions, we still have (\ref{func iden
proj2-d}) near $p^*$.
Now, for any $z'_0$ near 0, we can find a minimal rational curve
$C_{z'_0}$ passing through $\Psi_1^{-1}(z'_0, 0)$ such that the
tangent of $\Psi_1(C_{z'_0})$ at $(z'_0, 0)$ is transversal to $E$
at 0. (This is due to the fact that for any $q\in M$, the
holomorphic (isometric) isotropic subgroup acts irreducibly on
$T_q^{(1,0)}$ and the subspace of the span of the tangent vectors of
minimal rational curves at $q$ is an invariant subspace of
$T_q^{(1,0)}$.) Write a holomorphic parametrization of
$\Psi_l(C_{z'_0})$ as: $z'= \phi'(\eta), z_n= \phi_n(\eta)$ with $
\phi'(0)=z'_0, \phi_n(0) =0$ and  $\phi_n(\eta) = \eta h(\eta)$ with
$h(0)\not =0$. Restrict $Y^{*}$ to $\Psi_1(C_{z'_0})$ near $(z''_0,
z''_n)$ with certain $z''_n (\not= 0)\approx 0$ and $z''_0\approx
z_0'$. Let $\omega_{C_{z'_0}} = \omega_M \big|_{C_{z'_0}}$. Then we
know that $(C_{z'_0}, \omega_{C_{{z'_0}}})$ is isometric to
$(\mathbb{P}^1, \omega_1)$. Now, by Theorem \ref{kahler} with
$M={\mathbb P}^1$, we see that $F, G$ admit a global holomorphic
extension to $C_{z'_0}$. Hence,

$$\tilde{a}_{\alpha_0}(\phi'(\eta))h(\eta)^{\frac{\alpha_0}{N}} \equiv 0 ~~~\text{for}~~~\eta~\text{near}~0.$$
In particular, $\tilde{a}_{\alpha_0}(z'_0)=0$. Since $z'_0$ is
arbitrary, we see that $\tilde{a}_{\alpha_0}(z') \equiv 0$. This is
a contradiction.

\medskip

We next claim that $F$ (resp. $G$) admits a holomorphic extension
along any path $\gamma\subset M$ with $\gamma(0) \in U$. Indeed, if
not, we can find $t_0\in(0, 1]$ such that $F$ (resp. $G$) admits a
holomorphic extension along $\gamma \big|_{[0, t']}$ for any
$t'<t_0$ but not along $\gamma \big|_{[0, t_0]}$. Assume that
$\gamma(t_0) \in U_{i_0}$ for some $i_0$, where
$(U_{i_0},\Psi_{i_0})$ is a Nash algebraic  holomorphic chart of $M$
with the similar  property as described for $(U_1,\Psi_1)$. Then
$\gamma(t') \in U_{i_0}$ for $t'(<t_0)$ sufficiently close to $t_0$.
Composing with a holomorphic isometry of ${\mathbb P}^{n_0}$, we can
make $p^*:=\gamma(t')=[1,0,\cdots,0]$. Also, when $t'$ is
sufficiently close to $t_0$, we can assume that a small
neighborhood, still denoted by $U_{i_0}$,  of $\gamma(t')$ contains
$\gamma(t_0)$. Also, we still have a Nash holomorphic coordinate map
$\Psi_{i_0}$ over $U_{i_0}$. Notice that composing $F$ or $G$ on the
right by automorphisms of the target manifold will never change the
holomorphic extendability. Now composing $F_l$ and $G_j$ by
holomorphic isometres, we can  assume that $F_l(p^*)=[1,0,\cdots,
0]$ and $G_j(p^*)=[1,0,\cdots, 0]$. Hence, repeating exactly the
argument as above,
it follows that the map obtained by  restricting  the extended map
$F$ (resp. $G$) to a small neighborhood of $\gamma(t')$ admits a
holomorphic extension along curves  in $U_{i_0}$. In particular, we
conclude that $F$ and $G$ extend holomorphically along $\gamma([0,
t_0])$.
 This is a contradiction. Now,  since
$M$ is simply connected, we conclude that $F$ (resp. $G$) extends to
a holomorphic map $\tilde{F}$ (resp. $\tilde G$) from $M$ into
$\mathbb{P}^{N_1} \times \cdots \times \mathbb{P}^{N_m}$ (resp.
$\mathbb{P}^{N'_1} \times \cdots \times \mathbb{P}^{N'_v}$).

Finally, we show that each $\tilde{F}_l$  (
resp. $\tilde{G}_j$) is a holomorphic isometric embedding from $M$
into $\mathbb{P}^{N_l}$ (resp. $\mathbb{P}^{N'_j}$) up to a certain
isometric constant $m_l$ (resp. $n_j$) $\in
{\mathbb N}$.
Theorem \ref{kahler}, $\tilde{F}_l$,
  when restricted to each minimal rational curve $C$, is a
 holomorphic isometric embedding up to an isometric constant.
   Notice that $\tilde{F}_l ^* (\frac{1}{2\pi}\omega_{N_l})$ is a closed $(1, 1)$-form on $M$, which is also an element in $H^2(M_l,{\mathbb Z})$. Since $H^{2}(M, \mathbb{Z})$
   is generated by ${1\over 2\pi}\omega_M$,
    we have a certain $m_l\in {\mathbb Z}$ such that

$$\tilde{F}_l ^* \omega_{N_l} - m_l \omega_M = \sqrt{-1} \partial \bar\partial \vartheta$$ for a certain real-valued, real analytic function $\vartheta$
over $M$. Restricting to each minimal rational curve $C$, we get

$$\sqrt{-1}\partial\bar\partial \vartheta\big|_C = m_C \omega_C. $$ Since $\omega_C$ is positive definite,
we see that $\vartheta \big|_C$ is either subharmonic or
superharmonic over $C$. Thus $\vartheta$ is constant on $C$ as $C$
is compact. Since any two points in $M$ can be connected by a finite
sequence of minimal rational curves (see \cite{HK}), we see that
$\vartheta$ is a constant over $M$. Hence we have
$$\tilde{F}_l^* \omega_{N_l} = m_l \omega_M.$$
 Since $F_l$ is not constant, we see that $m_l>0$.
Now, applying the Nakagawa-Takagi theorem and applying the local
uniqueness theorem (up to isometries) of Calabi ([Theorem 9, C]), we
see that $F_l$ (resp. $G_j$) coincides with the (one-to-one)
$m_l$-th (resp. $n_j$-th) cannonical embedding from $M$ into $M_l$
(resp. $M_j'$) upto a unitary action. The
identity (\ref{james-0008}) thus also holds.
The proof of Theorem \ref{maintheorem-local} 
is complete. $\endpf$

\bigskip

\section{Examples and Remarks}

\bigskip \noindent {\bf Remark 6.1} (Calabi [p.23, C]):
  $({\mathbb P}^n, \mu\omega_n)$
can be locally holomorphically and  isometrically embedded into
$({\mathbb P}^\infty, \lambda\omega_\infty)$ if and only if
$\lambda=k\mu$ with $k\in {\mathbb N}$. This fact can  be easily
seen to be equivalent to the following simple  statement: {\it There
is a sequence of holomorphic functions $\{f_i(z)\}_{i=1}^{\infty}$
defined in a certain fixed small neighborhood of $0\in {\mathbb
C}^n$ with $f_i(0)=0$ for each $i$ such that
$(1+|z|^2)^{\lambda/\mu}=1+\sum_{i=1}^{\infty}|f_i(z)|^2$ if and
only if $\lambda/\mu$ is a positive integer.} (When $\lambda/\mu$ is
not a positive integer,  many coefficients for  terms of the form
$|z|^{2k}$ in the Taylor expansion of the left hand side are
negative, while such  kind of terms in the right hand side have
non-negative coefficients.)  Hence, as in  Calabi's paper  [p.23,
C], if $\mu_1/\mu_2$ is not rational, there is no $\lambda$ such
that both  $({\mathbb P}^n, \mu_1\omega_n)$ and $({\mathbb P}^m,
\mu_2\omega_m)$ can be (locally holomorphically and isometrically)
embedded into $({\mathbb P}^\infty, \lambda\omega_\infty)$. Thus
$({\mathbb P}^n\times {\mathbb P}^m, \mu_1\omega_n \oplus
\mu_2\omega_m)$ can not be locally holomorphically and isometrically
embedded into $({\mathbb P}^\infty, \lambda\omega_\infty)$ for any
choice of $\lambda$.

\bigskip

{\bf Example 6.2}: Let $\{ \mu_1, \cdots, \mu_m \}$ and
$\{\lambda_1, \cdots, \lambda_v\}$ be two sets of positive numbers.
Suppose that there exist nonnegative integers $m'_l, n'_j$ such that
$$\sum_{l=1}^m m'_l \mu_l = \sum_{j=1}^v n'_j \lambda_j
>0.$$ Now, if there is a  holomorphic map $G=(G_1,\cdots, G_v):
({\mathbb C}^n\subset {\mathbb P}^n,\omega_n)\ra ({\mathbb
P}^{N'_1}\times\cdots\times {\mathbb P}^{N'_v},\oplus_{j=1}^v
\lambda_j \omega_{N'_j})$
 and
 holomorphic map
$F=(F_1,\cdots, F_m): ({\mathbb C}^n\subset {\mathbb
P}^n,\omega_n)\ra ({\mathbb P}^{N_1}\times\cdots\times {\mathbb
P}^{N_m},\oplus_{l=1}^m \mu_l \omega_{N_l})$ such that
$\omega_{{\mathbb P}^n,G,\lambda}=\sum_{l=1}^m\mu_l
F^*_l\omega_{N_l},$  and each mapping factor is itself an isometry
up to a conformal factor, then we  see  that there are
 positive integers $m_l, n_j$ such that
\begin{equation}\label{james-11}
 \sum_{l=1}^m m_l
\mu_l = \sum_{j=1}^v n_j \lambda_j +1.
\end{equation}

On the other hand, assume (\ref{james-11}).  Given any holomorphic
map $f: 0\in U\subset\mathbb{C}^n \rightarrow \mathbb{C}^{N}$ with
$f(0)=0$, define

$$h_l=(1+|z|^2)^{m_l} (1+ |f|^2)^{m'_l}  ~~\text{and}~~q_j=(1+|z|^2)^{n_j} (1+ |f|^2)^{n'_j}.$$
Then
\begin{equation}\label{eg}
h_1^{\mu_1} \cdots h_m^{\mu_m} q_1^{-\lambda_1} \cdots
q_v^{-\lambda_v} = 1+ |z|^2.
\end{equation}
By  Lemma \ref{new-lemma}, one can construct holomorphic maps $F_l:
U\subset{\mathbb C}^n \subset \mathbb{P}^n \rightarrow
\mathbb{P}^{N_l}$ and $G_j: U\subset {\mathbb C}^n \subset
\mathbb{P}^n \rightarrow \mathbb{P}^{N'_j}$ such that
$(1+|F_l|^2)=(1+|z|^2)^{m_l} (1+ |f|^2)^{m'_l}$,
$(1+|G_j|^2)=(1+|z|^2)^{n_j} (1+ |f|^2)^{n'_j}$. Hence, by
(\ref{eg}), we have

$$\sum_{l=1}^m \mu_l F^*\omega_{N_l} = \sum_{j=1}^v \lambda_j G^*\omega_{N'_j} + \omega_n.$$
Notice that $f$ will be merged as part of the components of $F$ and
$G$. Since $f$ is arbitrarily assigned, we do not have algebraicity,
global extendability and rigidity for $F$ and $G$.

On the other hand, it is easy to see that when $\{\mu_l\}$
$\{\lambda_j\}$ satisfy (\ref{james-10}) and (\ref{james-0008}), for
each Hermitian symmetric space of compact type $(M,\omega_M)$, we
can find $F_l: (M,\omega_m)\ra ({\mathbb P}^{N_l},\omega_{N_l})$ and
$G_j: (M,\omega_m)\ra ({\mathbb P}^{N_j'},\omega_{N_j'})$ such that
$F^*_l\omega_{N_l}=m_l\omega_M$, $G^*_j\omega_{N_j'}=n_j\omega_M$,
and $\omega_{M,G,\lambda}=\sum_{l=1}^{m}\mu_l F_l^*\omega_{M_l}$.

 Finally, one can easily construct many examples of
$\lambda_j$ and $\mu_l$ such that both  (\ref{james-10}) and
(\ref{james-0008}) hold. One simple example is given as follows: Let
$v=1,m=2$, $\lambda_1 = \sqrt{2}, (\mu_1,
\mu_2)=(\sqrt{2}+\frac{1}{4}, \frac{1}{4})$. Then (\ref{james-10})
holds trivially. Meanwhile, $2\mu_1+2\mu_2=2\lambda_1+1.$ Also,
notice that $\mu_2/\mu_1$ is irrational and thus $({\mathbb
P}^{N_1}\times{\mathbb
P}^{N_2},\mu_1\omega_{N_1}\oplus\mu_2\omega_{N_2})$ can not be
embedded into $({\mathbb P}^{\infty},\mu\omega_{\infty})$ for any
$\mu>0.$

\bigskip

{\bf Example 6.3}:
Let $\{\mu_l\}_{l=1}^{m}, \{\lambda_j\}_{j=1}^v$ be two sets of positive
real numbers and let $\{m_l\}_{l=1}^{m}, \{n_j\}_{j=1}^v$ be two sets of positive nature
numbers such that (\ref{james-0008}) holds. Then for any irreducible
Hermitian symmetric space $(M,\omega_M)$ of compact type, equipped
with  a K\"ahler-Einstein metric $\omega_M$ normalized as in Theorem \ref{maintheorem-local}, by the Nakagawa-Takagi theorem [M1], there is a holomorphic
isometric embedding $F_l$ (resp. $G_j$) from $(M,\omega_M)$ into $({\mathbb
P}^{N_l},\omega_{N_l})$ (resp. (${\mathbb P}^{N'_j},\omega_{N'_j})$) with $\omega_{N_l}$ (resp. $\omega_{N'_j}$) the standard
Fubini-Study metric such that $F_l^{*} \omega_{N_l}=m_l\omega_M$ (resp. $G_j^{*} \omega_{N'_j}=n_j\omega_M$).
Thus, $F=(F_1,\cdots, F_m): (M,\omega_M+ \sum_{j=1}^v \lambda_j G^*_j \omega_{N'_j})\rightarrow ({\mathbb
P}^{N_1}\times \cdots\times {\mathbb P}^{N_m}, \oplus_{l=1}^m\mu_l
\omega_{N_l}) $ is an isometric embedding.


\noindent $^{*}$ Department of Mathematics, Rutgers University, New Brunswick, NJ 08903, USA\\

\noindent $^\dagger$ Department of Mathematics, Johns Hopkins University, Baltimore, MD 21218, USA\\

\end{document}